\documentclass[11pt]{article}
\usepackage[margin=1in]{geometry}
\usepackage[numbers,sort&compress]{natbib}

\usepackage{graphicx}\usepackage{multirow}\usepackage{amsmath,amssymb,amsfonts}\usepackage{amsthm}\usepackage{mathrsfs}\usepackage[title]{appendix}\usepackage{textcomp}\usepackage{booktabs}\usepackage{algorithm}\usepackage{algorithmicx}\usepackage{algpseudocode}\usepackage{enumitem}
\usepackage{bm} 
\usepackage{bbm}
\usepackage{url}
\usepackage[hidelinks]{hyperref}
\usepackage{caption}
\usepackage{subcaption}
\usepackage{array}
\usepackage{tabularx}

\usepackage{empheq}

      \newlength{\weightingTableWidth} 
\theoremstyle{plain}\newtheorem{theorem}{Theorem}\newtheorem{corollary}[theorem]{Corollary}
\newtheorem{lemma}[theorem]{Lemma}

\theoremstyle{remark}\newtheorem{remark}{Remark}
\theoremstyle{definition}

\DeclareMathAlphabet\EuScript{U}{eus}{b}{n}
\SetMathAlphabet\EuScript{bold}{U}{eus}{b}{n}

\def\mA{{\bm{A}}}

\def\mX{{\bm{X}}}

\def\vzero{{\bm{0}}}

\def\va{{\bm{a}}}
\def\vb{{\bm{b}}}

\def\vc{{\bm{c}}}

\def\ve{{\bm{e}}}

\def\vr{{\bm{r}}}

\def\vv{{\bm{v}}}

\def\vx{{\bm{x}}}

\def\vz{{\bm{z}}}

\author{Jessie Chen\thanks{Department of Mathematics, North Carolina State University, Raleigh, NC, USA. Email: jchen228@ncsu.edu.}\textsuperscript{,}\kern3pt\thanks{These authors contributed equally to this work.} \and Ziyu Cheng\thanks{Department of Computational Mathematics, Science and Engineering, Michigan State University, East Lansing, MI, USA. Email: chengzi3@msu.edu.}\textsuperscript{,}\kern3pt\footnotemark[2] \and Yunhui He\thanks{Department of Mathematics, University of Houston, Houston, TX, USA. Email: yhe43@central.uh.edu} \and Alona Kryshchenko\thanks{Department of Mathematics, Data Science and Statistics, California State University Channel Islands, Camarillo, CA, USA. Email: alona.kryshchenko@csuci.edu} \and Jocelyn T. Chi\thanks{Division of Computational Health Sciences, University of Minnesota Twin Cities, Minneapolis, MN, USA. Email: jocelync@umn.edu} \and Anna Ma\thanks{Department of Mathematics, University of California, Irvine, Irvine, CA, USA. Email: anna.ma@uci.edu} \and Jamie Haddock\thanks{Department of Mathematics, Harvey Mudd College, Claremont, CA, USA. Email: jhaddock@hmc.edu}}
\date{}

\begin{document}

\title{Weighted Randomized Kaczmarz for Sparsely Corrupted Linear Systems}

\maketitle

\begin{abstract}
In this paper, we introduce a generalized weighted randomized Kaczmarz framework to approximate solutions to overdetermined linear systems whose measurements are perturbed by sparse, additive corruptions. The proposed framework samples equations with probabilities proportional to a user-defined weighting function of their corresponding residual magnitudes. Our work focuses on weighting functions that will reduce the influence of potentially corrupted equations, and provides analyses for a general class of non-increasing weighting functions. We derive one-step expected squared error bounds that separate progress on uncorrupted equations from the contribution of corrupted equations, and provide sufficient conditions for geometric convergence under one specific family of weighting functions -- threshold weighting. We also investigate adaptive exponential weights scaled by the median or median absolute deviation of residual magnitudes, requiring no prior knowledge of the corruption fraction or magnitudes. Lastly, we consider a family of inverse residual weighting functions. Numerical experiments demonstrate the benefit of such a generalized framework in its flexibility in developing approaches for settings in which information about the corruptions is known or unknown.
\end{abstract}

\par

\par

\section{Introduction}\label{sec:intro}

Large-scale systems of linear equations arise in numerous data science, imaging, and network problems~\cite{GBH70:Algebraic-Reconstruction,HM93:Algebraic-Reconstruction,savvides2001dynamic}, frequently appearing as subroutines in machine learning~\cite{martens2010deep,green1984iteratively} and in several optimization methods~\cite{shi1995solving}.  In cases involving extremely large consistent or nearly-consistent systems, iterative methods are particularly attractive for their light computational and memory footprints and theoretical convergence guarantees. 
The Kaczmarz methods comprise a widely popular family of \emph{row-action} iterative methods, which consist of sequential orthogonal projections towards the solution set of a single equation~\cite{Kac37:Angenaeherte-Aufloesung}.
These methods gathered much renewed interest after the elegant convergence analysis of the \emph{Randomized Kaczmarz (RK) method} in~\cite{SV09:Randomized-Kaczmarz}. There, the authors showed that, under specific sampling distributions, the RK converges at least linearly in expectation for consistent systems with a unique solution.

Despite the rapidly increasing size of linear systems in modern applications, data corruptions that are difficult to detect and that can severely compromise accuracy remain a persistent problem. Examples of such corruptions include sensor malfunctions, transmission errors, storage contamination, or even adversarial perturbations. Standard iterative methods for linear systems are generally not robust to data corruptions. Rather, their focus is the \emph{noisy} case, where perturbations to the data are small in magnitude and uniformly spread across the system~\cite{needell2010randomized, zouzias2013randomized}.  

In this work, we focus instead on \emph{corrupted} linear systems of equations. Here, we assume that a small fraction of equations contain additive corruptions that can be large in magnitude or may even be adversarially designed. Suppose we are given an observed matrix $\mA \in \mathbb{R}^{m \times n}$ and right-hand side $\vb \in \mathbb{R}^{m}$ with $m \gg n$. We consider solving large systems of equations defined by $\mA$ and $\vb$, where $\vb = \mA \vx^* + \vc$ is the sum of a consistent right-hand side $\mA \vx^*$ and an additive sparse corruption vector $\vc$. The goal is to approximate the solution $\vx^*$ to the unobserved consistent system in the presence of a sparse corruption vector $\vc$. 

Recent works on iterative methods for solving corrupted linear systems include the Quantile Randomized Kaczmarz (QRK) method~\cite{haddock2022quantile} and its variants~\cite{cheng2023block,haddock2023subsampled}. These 
techniques restrict the update in each iteration to equations whose corresponding absolute residuals are less than a sufficiently-small quantile statistic of the entire absolute residual vector. Therefore, they require careful pre-selection of the quantile statistic to reduce the probability of sampling corrupted equations. Accordingly, these methods are very sensitive to the estimated corruption proportion. 

This raises our motivating question: \emph{Can we design row-sampling rules that are robust to sparse corruption without requiring any a priori knowledge of the additive corruption vector $\vc$?} 
In this work, we introduce a novel framework for developing variants of randomized iterative methods that are robust to corruption in the linear system measurements. We achieve this by employing a function of the residuals to downweight observations with large residual entry magnitudes. This diminishes the probability that the method samples potentially corrupted equations in each iteration. Our framework includes families of methods with varying levels of sensitivity to parameter tuning, enables flexibility to adjust how aggressively the method avoids sampling potential corruptions, and illustrates a trade-off at play: more conservative methods converge more slowly. 

In the rest of this section, we begin by setting some notation that we will refer to throughout this work (Section~\ref{sec:notation}). We then highlight our contributions (Section~\ref{sec:contributions}) and briefly review some related works (Section~\ref{sec:relatedworks}). 

\subsection{Notation.}\label{sec:notation} Throughout the work, we use superscripts to denote iteration indices.  We represent matrices with boldface capital Roman letters (e.g., $\mX$), vectors with boldface lower-case Roman letters (e.g., $\vx$), and (integer or real) scalars as unbolded lower-case Roman and Greek letters (e.g., $x$ or $\lambda$). For the matrix $\mA$, we write its $i$th row as $\va_i^\top$, and for a vector $\vv$, we denote its $i$th entry by $v_i$. $\mA^\dagger$ denotes the Moore-Penrose pseudo‑inverse of matrix $\mA$.  We also let $$s_{\max}(\mA):= s_1(\mA) \ge s_2(\mA) \ge \cdots \ge s_{\min}(\mA) := s_r(\mA) > s_{r+1}(\mA) = \cdots = s_n(\mA) = 0$$ be the singular values of matrix $\mA \in \mathbb{R}^{m \times n}$ with rank $r \le n$. $\| \cdot\|$ denotes the $\ell_2$ norm for both vectors and matrices.  

We let $[m]$ denote the set of the first $m$ positive integers. 
We denote the corruption ratio as $\beta$, where $0 < \beta < 1$.
We denote the set of indices of corrupted equations (the set of nonzero entries of $\vc$) as $C = \text{supp}(\vc)$, and let $|C|=\lfloor \beta m \rfloor$ denote its number of elements, where $\lfloor \beta m \rfloor$ is the maximum integer that is less or equal to $\beta m$. For an approximation $\vx$ to $\vx^{\star}$ 
let $\ve: = \vx - \vx^*$ and $\vr:= \mA\vx-\vb$ denote the error and residual vectors, respectively. We refer to the entries of $\vr$ as the \emph{residuals}. While denote the residual entry corresponding to the $i$th row with $r_{i}$, we also employ $r$ to refer generically to any entry in $\vr$.

For $C \subset [m]$, we denote the complement $C^c := [m]\setminus C.$ 
For a vector $\vz \in \mathbb{R}^m$ with sorted entries
$z_{(1)} \le \cdots \le z_{(m)}$, its \emph{empirical $\tau$-quantile}
is $Q_{\tau}(\vz) = z_{(\lceil m\tau \rceil)}$, where $\tau \in (0,1)$ and for any $x\ \in \mathbb{R}$, $\lceil x \rceil$ denotes the smallest integer greater than or equal to $x$. Therefore, $Q_{\tau}(\vz)$ is the smallest entry in $\vz$ so that at least a proportion $\tau$ of the entries in $\vz$ do not exceed it. 
Note that $Q_{0.5}(\vz)$ denotes the median of the entries in $\vz$.  

For a vector $\vv \in \mathbb{R}^k$ and scalar $\lambda$, we define the notation $\vv - \lambda := \vv - \lambda \mathbf{1}$ where $\mathbf{1} \in \mathbb{R}^k$ is the all-ones vector of length $k$. The element-wise absolute value of a vector is denoted by $|\vv|$. We denote the \emph{median absolute deviations} from the median of a vector as $\text{MAD}(\vv) := Q_{0.5}(|\vv - Q_{0.5}(|\vv|)|)$,
and rescale assuming a normal distribution to define $$\sigma_{\text{MAD}}(\vv) := 1.4826 \cdot\text{MAD}(\vv).$$

\subsection{Contributions}\label{sec:contributions}

We summarize our main contributions below.
\begin{enumerate}
    \item \textbf{A unifying weighted framework.} We introduce the
    generalized weighted RK method (GWRK,
    Algorithm~\ref{alg:gwrk}), a variant of the RK method which samples rows with probability
    proportional to a user-specified function $f$ of the residual
    magnitudes. GWRK includes Motzkin's method~\cite{motzkinschoenberg,agmon}, greedy RK~\cite{bai2018greedy},
    weighted RK~\cite{steinerberger2021weighted}, and QRK~\cite{haddock2022quantile} as special cases. Unlike prior
    weighted schemes, which favor large residuals to accelerate
    convergence, we use the framework to design weights that
    \emph{avoid} rows corresponding to corrupted equations.

    \item \textbf{A one-step error bound for general non-increasing weights.} Theorem~\ref{thm:nonincreasing_f_main} bounds the expected squared error
    after one iteration of GWRK applied with any non-increasing $f$. The bound separates
    a contraction term from a convergence-horizon term governed by the
    weight $f$ places on rows corresponding to corrupted equations. We identify cases when a horizon is
    unavoidable (Remark~\ref{rem:horizon}) and two ways to avoid one:
    thresholding (Remark~\ref{rem:threshold}) and weights whose relative
    mass on corrupted rows vanishes (Remark~\ref{rem:exponential}).

    \item \textbf{Three weighting families with specialized guarantees.}
    We specialize Theorem~\ref{thm:nonincreasing_f_main} to three families of weighting functions in Subsection~\ref{subsec:special_cases}: thresholded 
    (Corollary~\ref{cor:thresholded_convergence}), exponential 
    (Corollary~\ref{cor:medianexp}), and inverse-magnitude weighting functions
    (Corollary~\ref{cor:inverse-mag-weight}). 
    \begin{itemize}
        \item For thresholded  weighting functions, we prove geometric convergence with
        no convergence horizon under explicit conditions. This result applies to QRK
        and variants with adaptive quantile levels.
        \item For exponential weighting  functions, we prove bounds when $f$ is scaled by 
        the median $Q_{0.5}$ or $\sigma_{\mathrm{MAD}}$ of the absolute residuals      (Lemma~\ref{lem:exponential_function_bounds}). Unlike QRK and GWRK with other threshold weighting functions, GWRK with these exponential weighting functions does not require estimating the corruption rate
        $\beta$ or magnitudes.
        \item For inverse-magnitude weighting functions, we provide a one-step bound
        that may explain the early convergence seen in practice for these
        highly conservative methods.
    \end{itemize}
    Together, these results describe a trade-off: the more conservative a
    weighting function is about sampling potentially corrupted rows, the less
    prior knowledge it needs about the corruption setting. However, the trade-off is that it becomes more challenging 
    to obtain a
    horizon-free guarantee. Indeed, we see that the most conservative weighting functions often suffer from extremely slow convergence. 

    \item \textbf{A link to robust statistics.} We make connections to robust regression works from the statistics literature (Section~\ref{subsubsec:robustreg}). In particular, much of the robust statistics literature relies on loss functions that comprise functions of residuals. Additionally, the exponential weights we introduce later are inspired by the weights employed in $L^2$ distance estimation. 

    \item \textbf{Numerical experiments.}     We test example functions from each of the three weighting families, with results reported in Section~\ref{sec:experiments}.
    Threshold weighting converges rapidly when an appropriate threshold is chosen.
    The parameter-free exponential functions also show good convergence behavior. Although they are more conservative and generally require more iterations to converge than QRK, some exhibit similar convergence slopes and evolution dynamics to the quantile threshold weighting.
    The inverse-magnitude weighting is the most conservative of the three families and exhibits extremely slow convergence. We run an additional numerical experiment on the exponential weighting family of functions by varying the corruption ratio $\beta$ to show that the restriction on $\beta$ in Corollary~\ref{cor:medianexp} is empirically necessary. That is, when $\beta > 0.5$ we begin to observe a convergence horizon and when $\beta < 0.5$ we observe convergence, so the assumption $\beta < 0.5$ is necessary for guarantees of error contraction. Surprisingly, the method converges more quickly with corruption ratio $\beta$ closer to $0.5$, since the algorithm becomes restricted to a smaller uncorrupted system, which offers faster convergence due to more frequent repeat sampling across these solution-defining rows.

\end{enumerate}

\subsection{Related Work}\label{sec:relatedworks}

In this section, we begin by reviewing the RK method, followed by the QRK method. We also highlight some works from the robust statistics literature that have particular connections to our work here. 

\subsubsection{Randomized Kaczmarz method} 
The Kaczmarz method~\cite{karczmarz1937angenaherte}, an alternating-projections method, has grown popular in recent years for solving large, consistent linear systems.
Given the system $\mA\vx = \vb$, the method computes iterates by projecting onto the hyperplane defined by the equation $\va_i^\top \vx = b_i$. The iterate in the $t$th iteration is recursively
defined as  $$\vx^{(t)} = \vx^{(t-1)} + \frac{b_{i_t} - \va_{i_{t}}^\top \vx^{(t-1)}}{\|\va_{i_{t}}\|^2} \va_{i_t}.$$ The index $i_t$ is usually either selected cyclically (i.e., $i_{t+1} = i_t + 1$) or randomly (i.e., $i_t \sim \mathcal{D}$ where $\mathcal{D}$ is a distribution over $[m]$).  Note that the method can be initialized with any vector $\vx^{(0)} \in \mathbb{R}^n$.

In~\cite{SV09:Randomized-Kaczmarz}, the authors showed that for a consistent system with the unique solution $\vx^*$, the RK method (with specified sampling distribution) converges at least linearly in expectation with the guarantee
\begin{equation}
    \mathbb{E}\|\ve^{(t)}\|^2 \le \left(1 - \frac{s_{\min}^2(\mA)}{\|\mA\|_F^2}\right)^t \|\ve^{(0)}\|^2, \label{eq:RKrate}
\end{equation} 
where $\ve^{(t)} := \vx^{(t)} - \vx^*$. 
 Subsequently, many variants of RK have been proposed; see e.g.,~\cite{bai2018greedy,bai2021greedy,needell2015randomized,necoara2019faster,niu2020greedy,moorman2021randomized,tondji2023faster,du2020randomized,cheng2023block,haddock2023subsampled,steinerberger2023quantile,zeng2023randomized,eldar2011acceleration,jiang2026accelerated}. The applications of RK-type methods have been extended to least squares~\cite{zouzias2013randomized}, tensor problems~\cite{wang2022randomized,castillo2025randomized,chen2021regularized}, and nonlinear systems~\cite{wang2022nonlinear,liu2025greedy}. 

One direction of particular interest in the literature is developing row sampling techniques that select an index that produces a large magnitude update and hopefully speeds convergence; these are often referred to as \emph{greedy} methods and utilize sampling rules that bias towards large magnitude residuals.   The greedy randomized Kaczmarz (GRK) method~\cite{bai2018greedy} balances between a greedy selection and a random selection by randomly selecting among a subset of greedily chosen rows. At every iteration, the residual is computed and then used to calculate the effective threshold for $\omega \in (0,1]$: 
$$ \epsilon_t = \omega \max_i \frac{|r^{(t)}_{i}|^2}{\|\vr^{(t)} \|^2}+ (1-\omega)  \frac{1}{m}, $$
assuming $A$ is row normalized. Then, in iteration $t$ the sampled row index $i_t \in [m]$ is selected with probability 
\begin{equation}\label{eq:defintion-pti}
p_{i} \propto\begin{cases} |r^{(t)}_i|^2 & \text{if } |r^{(t)}_i|^2/\|\vr^{(t)} \|^2 \ge \epsilon_t , \\ 
        0 & \text{otherwise.}\end{cases}
\end{equation}
The convergence result in~\cite{bai2018greedy} shows that the method is expected to converge, though it generally offers no significant benefit over RK.
\begin{theorem}(Main Theorem~\cite[Theorem 3.1]{bai2018greedy}) Let $A$ be normalized to have unit norm rows and $x^{(t)}$ be iterates of the GRK algorithm. Then
\begin{equation*}
    \mathbb{E} \|\vx^{(1)} - \vx^* \|^2 \leq \Big(1 - \frac{s^2_{\min}(\mA)}{m} \Big) \| \vx^{(0)} - \vx^*\|^2,
\end{equation*}
and for $t \geq 1$
\begin{equation*}
    \mathbb{E} \|\vx^{(t+1)} - \vx^* \|^2 \leq \left(1 - \frac{1}{2} \left(\frac{m}{m-1} 
    +1\right) \frac{s^2_{\min}(\mA)}{m}\right) \| \vx^{(t)} - \vx^*\|^2.
\end{equation*}
\end{theorem}

If $\omega = 1$, GRK corresponds to Motzkin's method, which only chooses the row corresponding to the largest residual~\cite{motzkinschoenberg,HN18Motzkin}. 
The maximal weighted residual randomized Kaczmarz method proposed in~\cite{du2019new} is a variant of Motzkin's method which instead selects the row with largest \emph{relative} residual magnitude, $|\va_j^\top \vx^{(t-1)} - b_j|/\|\va_j\|$.

In~\cite{steinerberger2021weighted}, a weighted RK algorithm (referred to here as WRK) was introduced where in each iteration, rows of $\mA$ are sampled with probabilities
$$p_{i} = \frac{| \langle \va_i, \vx^{(t)} \rangle - b_i |^q}{\|\mA \vx^{(t)} - \vb \|^q_{\ell_q}},$$ where $0 < q < \infty$ is a user-chosen parameter. The following result shows that the convergence will be \textit{at least} the rate of RK.   
\begin{theorem}[Main Theorem~\cite{steinerberger2021weighted}] Let $0 < q < \infty$ and $\mA$ be normalized to have unit norm rows. Then the iterates of WRK satisfy
\begin{equation*}
    \mathbb{E} \|\vx^{(t)} - \vx^* \|^2 \leq \left(1 - \inf_{\vz \neq 0} \frac{\| \mA\vz\|_{\ell_{q+2}}^{q+2}}{\| \mA\vz\|_{\ell_{q}}^{q}\| \| \vz\|^2}\right)^t \|\vx^{(0)} - \vx^* \|^2,
\end{equation*}
    and 
$$ \inf_{\vz \neq 0} \frac{\| \mA\vz\|_{\ell_{q+2}}^{q+2}}{\| \mA\vz\|_{\ell_{q}}^{q}\| \| \vz\|^2} > \frac{1}{\| \mA\|^2_F \| \mA^\dagger \|^2}.$$
\end{theorem}
In particular, as $q \rightarrow \infty$, this weighting scheme approximates Motzkin's method and only selects the row corresponding to the largest absolute residual. 

In this work, we consider a weighted RK for solving sparsely corrupted linear systems. Unlike previous works that consider the weighted RK, our work presents a general framework for reweighting, which provides flexibility to different settings, including robustness to sparse corruption, and includes GRK and WRK as special cases.  Our setting also differs from the primary setting of these methods in two ways.  
First, we do not seek the least squares solution. Rather, our goal is to compute a solution that fully avoids the effect of additive sparse corruption in the measurements. Second, our goal for our weighting framework is not to speed convergence, but only to avoid the effects of sparse corruption.  Greedy and random-greedy methods tend to assign weights proportional to a function of the residual, which can sometimes lead to an acceleration in convergence. However, when systems are sparsely corrupted, large residuals can correspond to outliers, rendering these greedy methods ineffective.

\subsubsection{Quantile randomized Kaczmarz methods}\label{subsec:QRK}

In many practical applications, the underlying large-scale linear systems are inconsistent because of sparse corruptions in the measured right-hand side. In such cases, the classical RK method becomes ineffective or may even fail to converge. One remedy is the QRK method~\cite{haddock2022quantile}, which is specifically designed to mitigate the influence of corrupted measurements. The key idea behind QRK is to perform the usual Kaczmarz projection step only with rows whose corresponding residual magnitude falls below a chosen quantile of all residual magnitudes, indicating that the sampled equation is unlikely to be corrupted. 
Since its introduction, improved analyses of QRK have been provided~\cite{steinerberger2023quantile} and several variations of the QRK method have been proposed. These include block accelerations~\cite{cheng2023block, zhang2025quantile}, lower quantile accelerations~\cite {battaglia2025reverse}, subsampled methods~\cite{cai2025subsample,haddock2023subsampled, battaglia2026quantile}, and many others~\cite{zhang2025quantile, coria2024quantile}. Algorithm~\ref{alg:QRK} presents one implementation of the QRK method~\cite{steinerberger2023quantile}. 
\begin{algorithm}[H]
\caption{Quantile RK (QRK)} \label{alg:QRK}
\begin{algorithmic}
\State {\bf Inputs:} $\mA, \vb$, quantile $\tau$, maximum number of iterations $T$ 
\State $\vx^{(0)} = \vzero$
\For{$t = 1, \dots, T$}
    \State Sample $i_t \sim \text{Unif}\left(\left\{k \in [m]: |\langle \va_{k}, \vx^{(t-1)}\rangle -b_{k}| \le Q_{\tau}
    \left( | \mA \vx^{(t-1)} - \vb | \right)\right\}\right)$
        \State $\vx^{(t)} = \vx^{(t-1)} - \frac{\langle \va_{i_t}, \vx^{(t-1)}\rangle - b_{i_t} }{\|\va_{i_t}\|^2}\va_{i_t}$
\EndFor 
\State {\bf Return} $\vx^{(T)}$
\end{algorithmic}
\end{algorithm}

\subsubsection{Robust regression techniques}
\label{subsubsec:robustreg}

Although our work has roots in the RK method, we also draw connections to the literature on robust regression in statistics. Least squares regression estimates the mean of the response, or right-hand side $\vb$, conditional on the observed data $\mA$. Since mean estimates are sensitive to outliers, much of the robust regression literature focuses on identifying outlying observations, or rows, and reducing their influence on the solution.

The conventional approach is to diminish the weights of observations with large absolute residuals. Examples of these regression methods include least absolute deviations (LAD)~\cite{bloomfield1983least, Bassett01091978}, which minimizes the sum of the absolute residuals in lieu of their squares; least median of squares~\cite{rousseeuw1984least, stromberg1993computing}, which minimizes the median of the squared residuals; trimmed least squares~\cite{ruppert1980trimmed, agullo2001new}, which minimizes the sum of a subset comprising the smallest squared residuals; and trimmed LAD~\cite{honore1992trimmed}, which minimizes the sum of a subset of the smallest absolute residuals. These now have many variations and algorithmic adaptations. One example is an iterative trimmed least squares algorithm~\cite{shen2019iterative} for additive corruptions in the response. 

While the methods above all minimize sums of functions of the residuals or their trimmed variants, we highlight two lines of work from the robust statistics literature that are particularly connected to this work. The first is robust M-estimation \cite{susanti2014m, koenker1990m, huber1973robust}, which minimizes the sum of a function of the residuals $\rho(r_{i})$, where $r_{i}$ denotes the residual for the $i$th row. Choices for $\rho$ that are common for robust regression include the absolute loss $\rho(r_{i}) = |r_{i}|$ for LAD described previously, the Huber loss~\cite{huber1973robust, koenker1990m} $\rho_{\mathrm{H},c}$ and Tukey's biweight~\cite{susanti2014m}, or bisquare, loss $\rho_{\mathrm{T},c}$. The latter two are given by
\begin{eqnarray*}
\rho_{\mathrm{H},c}(r_{i}) &= \begin{cases} \frac{1}{2}r_{i}^2, & |r_{i}| \leq c, \\ c|r_{i}| - \frac{1}{2}c^2, & |r_{i}| > c.\end{cases}  \quad \text{ and} \quad
\rho_{\mathrm{T},c}(r_{i}) &= \begin{cases} r_{i}\left(1-(r_{i}/c)^2\right)^2, & |r_{i}| < c, \\ 0, & |r_{i}| \geq c. \end{cases}
\end{eqnarray*}
respectively, where $c>0$ is a threshold parameter. The Huber loss is quadratic for small residuals but resembles the absolute loss for larger ones. Meanwhile, Tukey's biweight is approximately quadratic for small residuals but is zero for larger ones.

A second related line of work that is particularly connected to this work is $L^2$ distance estimation (L$_2$E)~\cite{chi2022user, liu2023sharper, scott2009l2e}. Unlike classical likelihood-based methods, which can be highly sensitive to outlying rows, the L$_2$E employs the $L^2$ distance between the model and the observed data. This yields estimating equations in which each row has a weight $w_i = \exp \left(-\frac{r_{i}^{2}}{2\sigma}\right)$, where $\sigma$ is an estimate for the population standard deviation. These weights decrease with the magnitude of the residual so that observations with small residuals retain near-full influence while observations with large residuals are smoothly downweighted~\cite{chi2022user}. A key feature of the L$_2$E framework is that one can estimate $\sigma$ from the data and update the estimate in each iteration. 

Motivated by functions of residuals in M-estimation and smooth downweighting of outlying observations in the L$_{2}$E, we integrate these features into the weighted framework in WRK. In the next section, we introduce a framework for generalized weighted RK for corrupt linear systems with additive errors.

\section{A weighted framework for corrupt linear systems}

Inspired by the weighted framework in WRK~\cite{steinerberger2021weighted} and a sampling scheme to avoid corruptions in QRK~\cite{haddock2022quantile,steinerberger2023quantile}, our \textbf{Generalized Weighted Randomized Kaczmarz (GWRK)} method samples rows with probability proportional to a function $f$ of the absolute residuals.  Algorithm~\ref{alg:gwrk} shows pseudocode for the GWRK method.

\begin{algorithm}[H]
\caption{Generalized Weighted Randomized Kaczmarz (GWRK)} \label{alg:gwrk}
\begin{algorithmic}
\State {\bf Inputs:} $\mA, \vb, \vx^{(0)}, f(r;\vr), $ maximum number of iterations $T$
\For{$t = 1, \cdots, T$}
    \State{$\vr^{(t-1)} = \mA \vx^{(t-1)} - \vb$}
    \State{Sample row $i_t \in [m]$ with probability $p_{i_t} \propto f(|r^{(t-1)}_{i_t}|; |\vr^{(t-1)}|)$}.
    \State{$\vx^{(t)} = \vx^{(t-1)} - \frac{r^{(t-1)}_{i_t}}{\|\va_{i_t}\|^2} \va_{i_t}$}  
\EndFor
\State {\bf Return} $\vx^{(T)}$ 
\end{algorithmic}
\end{algorithm}

Algorithm~\ref{alg:gwrk} encompasses several existing methods as special cases.  For example, when the function $f$ assigns the largest magnitude residual probability 1 and all others zero, $$f(r; \vr) = \begin{cases} 1 & \text{ if } |r| = \|\vr\|_\infty \\ 0 & \text{ otherwise}\end{cases},$$ then we recover Motzkin's method~\cite{motzkinschoenberg,agmon,HN18Motzkin}.   Similarly, \eqref{eq:defintion-pti} provides exactly the function $f$ that recovers GRK~\cite{bai2018greedy} from GWRK.  
In the absence of data corruption, one could accelerate convergence by designing $f$ to favor sampling rows with larger absolute residuals.  Since we are working in the corrupt data regime, however, we seek to reduce the influence of contaminated rows. For this reason, we focus on non-increasing functions that assign larger weights to rows with smaller residual magnitudes.

We next prove our main theorem, which provides a bound on the expected squared error of GWRK after one iteration for a general non-increasing weighting function $f$.  This bound illustrates how the weighting function plays a role in both contraction of the error as well as a potential convergence horizon.  

For notational simplicity, throughout the remainder of this work, the weighting function $f(r; \vr)$ will be parameterized by the vector $|\vr^{(t-1)}|$ in the $t$th iteration, and so we introduce the notation $f_t(r) := f(r; |\vr^{(t-1)}|)$ for simplicity.  Here the variable $r$ represents an arbitrary scalar input to the function $f_t$; we apply $f_t$ always only to values that are nonnegative, so the behavior of $f_t$ on $r < 0$ is not relevant to our analysis.  Also, to ensure that Algorithm~\ref{alg:gwrk} is well-defined, we assume throughout that $f_t(|r^{(t-1)}_j|) > 0$ for at least one $j \in [m]$. We denote the conditional expectation conditioned on the outcome of the $(t-1)$th iterate with $\mathbb{E}_{t-1}[\cdot] := \mathbb{E}\left[\cdot \mid \vx^{(t-1)}\right]$.  

\begin{theorem}
    Suppose $\vb = \mA\vx^* + \vc$, where full-rank $\mA \in \mathbb{R}^{m \times n}$ with $m \ge n$ is row-normalized and $C = \text{supp}(\vc)$ satisfies $|C| = \beta m$ with $0 < \beta < 1$.     Suppose that the weighting function $f_t(r)$ is non-increasing and nonnegative on nonnegative inputs with $f_t(0)> 0$.  Let $\ve^{(t)} = \vx^{(t)} - \vx^*$ be the error in the $t^{th}$ iteration of Algorithm~\ref{alg:gwrk} applied to $\mA$ and $\vb$.  Define $c_{\min}^{(t-1)} = \max\{\min_{j \in C} |c_j|, \|\ve^{(t-1)}\|\}$ and assume $f_t(\|\ve^{(t-1)}\|) > 0$.  Then the errors satisfy
    \begin{align}
    \begin{split}
        \mathbb{E}_{t-1}[\|\ve^{(t)}\|^2] 
        &\le \left(1 - \frac{f_t(\|\ve^{(t-1)}\|) s_{\min}^2(\mA_{C^c})}{m f_t(0)}\right) \|\ve^{(t-1)}\|^2  \\
        &\quad\quad+\frac{f_t(c_{\min}^{(t-1)} - \|\ve^{(t-1)}\|)}{(1 - \beta)mf_t(\|\ve^{(t-1)}\|)} \|\vc\|^2.
    \end{split}
     \label{eq:main-bound}
    \end{align} \label{thm:nonincreasing_f_main} 
    
\end{theorem}

\begin{remark}
    We note that the constant on the convergence horizon term, i.e., the second term on the right-hand side of \eqref{eq:main-bound}, is lower-bounded with
    $$\frac{f_t(c_{\min}^{(t-1)} - \|\ve^{(t-1)}\|)}{(1-\beta)m f_t(\|\ve^{(t-1)}\|)} \ge \frac{f_t(c_{\min}^{(t-1)})}{(1-\beta)m f_t(0)}.$$
    The ratio $f_t(c_{\min}^{(t-1)})/f_t(0)$ reveals the best-case asymptotic likelihood of selecting the most difficult-to-detect corrupted row.  Recall that $\vr^{(t-1)} = \mA\vx^{(t-1)}-\vb=\mA \ve^{(t-1)} - \vc$.  Thus, if the iterates converge to $\vx^*$ and so the errors converge to $\mathbf{0}$, the corresponding residuals converge to $-\vc$.  In the distribution defined by $f_t$ applied to $|\vr^{(t-1)}|$, the uncorrupted rows have weights converging to $f_t(0)$, the maximum possible value of the weighting function, while the corrupted rows have weights converging to $f_t(|\vc_j|) \le f_t(c_{\min}^{(t-1)})$.  As long as $f_t(c_{\min}^{(t-1)}) > \epsilon > 0$, there is a fixed nonzero probability that a corruption may be selected, and thus, we cannot avoid a convergence horizon in this analysis for general choices of $f_t$.  \label{rem:horizon}
\end{remark}

\begin{remark}
    One way to avoid a convergence horizon in the analysis is to guarantee that there exists an iteration $T$ after which  $f_t(c_{\min}^{(t-1)} - \|\ve^{(t-1)}\|) = 0$ for all iterations $t \ge T$.  If such an iteration $T$ exists and $f_t(\|\ve^{(t-1)}\|)/f_t(0) \ge \delta > 0$ for all $t$, then the expected error of Algorithm~\ref{alg:gwrk} goes to 0 as $t \rightarrow \infty$.  A natural type of weighting function for the goal of guaranteeing that $f_t(c_{\min}^{(t-1)} - \|\ve^{(t-1)}\|) = 0$ is a \emph{thresholded} weighting function where there exists some $\theta^{(t-1)} > 0$ so that for any $r > \theta^{(t-1)}$, $f_t(r) = 0$.  The behavior of Algorithm~\ref{alg:gwrk} with this type of function is explored in Subsection~\ref{subsubsec:threshold_type}.  \label{rem:threshold}
\end{remark}

\begin{remark}
    Another way one might avoid a convergence horizon in this analysis is to utilize a sequence of weighting functions $f_t$ such that $f_t(c_{\min}^{(t-1)} - \|\ve^{(t-1)}\|)/f_t(\|\ve^{(t-1)}\|) \rightarrow 0$ as $t \rightarrow \infty$. In this case, since $f_t$ is non-increasing on inputs, this implies that $f_t(|r^{(t-1)}_i|) \rightarrow 0$ for all $i \in C$ as $t \rightarrow \infty$ and the probability of sampling any corrupted row goes to 0.  If the probability of sampling any corrupted row goes to 0 quickly enough, one would expect that any accrual of error through the second term in the bound~\eqref{eq:main-bound} in early iterations will be removed in later iterations.  We identify a promising class of \emph{exponential} weighting functions and explore the behavior of Algorithm~\ref{alg:gwrk} utilized with these $f_t$ in Subsection~\ref{subsubsec:exponential_type}. \label{rem:exponential}
\end{remark}

We now provide a proof of Theorem~\ref{thm:nonincreasing_f_main}.  Throughout this proof, we highlight several equations from which we will continue proofs of several special case corollaries in Subsection~\ref{subsec:special_cases}.  In the course of the proof, and later in the paper, we refer to disjoint `fast' and `slow' subsets of $C$, denoted as $C_F^{(t)}$ and $C_S^{(t)}$ respectively, and are defined as 
\begin{equation}
    C_F^{(t)} := \{j \in C : |c_j| > |\langle \va_{j}, \ve^{(t-1)}\rangle| \text{ and } f_t(|r_j^{(t-1)}|) > 0\} \label{eq:fast_C}
\end{equation} and 
\begin{equation}
    C_S^{(t)} := \{j \in C : |c_j| \le |\langle \va_{j}, \ve^{(t-1)}\rangle| \text{ and } f_t(|r_j^{(t-1)}|) > 0\}.\label{eq:slow_C}
\end{equation}
Here, $C_F^{(t)}$ is the set of indices of corrupted equations where the uncorrupted residual magnitude is less than the magnitude of the corruption affecting the equation -- we say it is `fast' because the residual error in the corresponding uncorrupted equation is small relative to the corruption, so one may hope that this corruption is now detectable relative to the residual errors of the other uncorrupted equations.  The slow set $C_S^{(t)}$ is the set of indices of the other corrupted equations -- these we expect we may not yet be able to detect due to the slow convergence of the uncorrupted residual error of these equations.

\begin{proof}
    We have $$\|\vx^{(t)} - \vx^*\|^2 = \|\vx^{(t-1)} - \vx^*\|^2 - |\langle \va_{j_t}, \vx^{(t-1)} - \vx^*\rangle|^2 + c_{j_t}^2.$$

    Now, taking expectation conditional on the first $t-1$ iterations, we have

    \begin{align}
        \mathbb{E}_{t-1}[\|\ve^{(t)}\|^2] &= \|\ve^{(t-1)}\|^2 - \mathbb{E}_{t-1}[|\langle \va_{j_t}, \ve^{(t-1)}\rangle|^2] + \mathbb{E}_{t-1}[c_{j_t}^2] \nonumber
        \\&= \|\ve^{(t-1)}\|^2 - \frac{1}{\sum_{i=1}^m f_t(|r_i^{(t-1)}|)} \sum_{j=1}^m f_t(|\vr_j^{(t-1)}|) |\langle \va_{j}, \ve^{(t-1)}\rangle|^2 
        \nonumber\\
        &\quad\quad\quad\quad\quad\quad\quad\quad\quad+ \frac{1}{\sum_{i=1}^m f_t(|r_i^{(t-1)}|)} \sum_{j=1}^m f_t(|r_j^{(t-1)}|) c_{j}^2\nonumber
        \\\begin{split} 
        &=\|\ve^{(t-1)}\|^2 - \frac{1}{\sum_{i=1}^m f_t(|r_i^{(t-1)}|)} \sum_{j\not\in C} f_t(|r_j^{(t-1)}|) |\langle \va_{j}, \ve^{(t-1)}\rangle|^2 
    \\&\quad\quad\quad\quad\quad\quad\quad\quad\quad+ \frac{1}{\sum_{i=1}^m f_t(|r_i^{(t-1)}|)} \sum_{j\in C} f_t(|r_j^{(t-1)}|) [c_{j}^2 - |\langle \va_{j}, \ve^{(t-1)}\rangle|^2]. 
        \end{split} \label{eq:three_terms}
    \end{align}

We proceed by lower bounding the second term and upper bounding the last term of ~\eqref{eq:three_terms}. To that end, note that due to the non-increasing property of $f_t$, we have $$mf_t(0) \ge \sum_{i=1}^m f_t(|r_i^{(t-1)}|) \ge \sum_{i\not\in C} f_t(|r_i^{(t-1)}|).$$ 
Then, we have 
\begin{align}
   \frac{1}{\sum_{i=1}^m f_t(|r_i^{(t-1)}|)} & \sum_{j\not\in C} f_t(|r_j^{(t-1)}|) |\langle \va_{j}, \ve^{(t-1)}\rangle|^2 \nonumber \\
   & \geq \frac{1}{m f_t(0)} \sum_{j \not\in C}  |\langle \tilde{\va}_{j}, \ve^{(t-1)}\rangle|^2 \nonumber \\
   & = \frac{1}{m f_t(0)} \|\tilde{\mA}_{C^c}(\ve^{(t-1)})\|^2\nonumber \\
   & \geq \frac{ s_{\min}^2(\tilde{\mA}_{C^c})}{m f_t(0)} \|\ve^{(t-1)}\|^2, \label{eq:lb-term2}
\end{align}
where $\tilde{\mA}_{C^c}$ has $j$th row $\tilde{\va}_j := \sqrt{f_t(|r_j^{(t-1)}|)} \va_j$.

  Since the terms corresponding to the `slow' set $C_S^{(t)}$ in the summation in the last term of~\eqref{eq:three_terms} are non-positive, $$\sum_{j \in C_S^{(t)}} f_t(|r_j^{(t-1)}|) [c_{j}^2 - |\langle \va_{j}, \ve^{(t-1)}\rangle|^2] \le 0,$$ we can drop these terms when deriving an upper bound for the last term in~\eqref{eq:three_terms}. Thus, 
\begin{align}
    \frac{1}{\sum_{i=1}^m f_t(|r_i^{(t-1)}|)} &\sum_{j\in C} f_t(|r_j^{(t-1)}|) [c_{j}^2 - |\langle \va_{j}, \ve^{(t-1)}\rangle|^2] \nonumber \\ 
    & \leq \frac{1}{\sum_{j\not\in C} f_t(|\langle \va_{j}, \ve^{(t-1)}\rangle|)} \sum_{j \in C_F^{(t)}} f_t(|r_j^{(t-1)}|) [c_{j}^2 - |\langle \va_{j}, \ve^{(t-1)}\rangle|^2] \nonumber \\ 
    & \leq \frac{f_t(c_{\min}^{(t-1)} - \|\ve^{(t-1)}\|)}{\sum_{j\not\in C} f_t(|\langle \va_{j}, \ve^{(t-1)}\rangle|)} \sum_{j \in C_F^{(t)}}  [c_{j}^2 - |\langle \va_{j}, \ve^{(t-1)}\rangle|^2]     \nonumber \\
    & \leq \frac{f_t(c_{\min}^{(t-1)} - \|\ve^{(t-1)}\|)}{\sum_{j\not\in C} f_t(|\langle \va_{j}, \ve^{(t-1)}\rangle|)} \|\vc\|^2, \label{eq:coarse_c_bound}
\end{align}
where the second inequality follows from the reverse triangle inequality applied to $j \in C_F^{(t)}$,
\begin{equation*}
    |r_j^{(t-1)}| \geq |c_j| - |\langle \va_{j}, \ve^{(t-1)}\rangle|  \geq c_{\min}^{(t-1)} - \|\ve^{(t-1)}\|,
\end{equation*}
which follows from Cauchy-Schwarz and the fact that $\mA$ has unit-norm rows as well as the definition of $c_{min}^{(t-1)}=\max\{\min_{j \in C} |c_j|, \|\ve^{(t-1)}\|\}$, i.e.,
\begin{equation}
|\langle \va_{j}, \ve^{(t-1)}\rangle| \le \|\va_{j}\|\|\ve^{(t-1)}\| = \|\ve^{(t-1)}\|.\label{eq:CS_unitnorm}
\end{equation}
Applying the bounds~\eqref{eq:lb-term2} and~\eqref{eq:coarse_c_bound} to \eqref{eq:three_terms}, we have 
\begin{equation}
    \mathbb{E}_{t-1}[\|\ve^{(t)}\|^2] \le \left(1 - \frac{s_{\min}^2(\tilde{\mA}_{C^c})}{m f_t(0)}\right) \|\ve^{(t-1)}\|^2 +\frac{f_t(c_{\min}^{(t-1)} - \|\ve^{(t-1)}\|)}{\sum_{j\not\in C} f_t(|\langle \va_{j}, \ve^{(t-1)}\rangle|)} \| \vc\|^2.
    \label{eq:almost_bound}
\end{equation}

    Now, beginning from~\eqref{eq:lb-term2}, we note that 
    \begin{align*}
    \|\tilde{\mA}_{C^c}(\ve^{(t-1)})\|^2 &= \sum_{j\not\in C} f_t(|\langle \va_{j}, \ve^{(t-1)}\rangle|) |\langle \va_{j}, \ve^{(t-1)}\rangle|^2 \\
    & \ge f_t(\|\ve^{(t-1)}\|) \sum_{j\not\in C} |\langle \va_{j}, \ve^{(t-1)}\rangle|^2 \\
    & \ge f_t(\|\ve^{(t-1)}\|) s_{\min}^2(\mA_{C^c})\|\ve^{(t-1)}\|^2,
    \end{align*}
    where the first inequality follows from the nonincreasing property of $f_t$ and~\eqref{eq:CS_unitnorm}.  
    Thus, we may replace $s_{\min}^2(\tilde{\mA}_{C^c})$ by $f_t(\|\ve^{(t-1)}\|) s_{\min}^2(\mA_{C^c})$ in~\eqref{eq:almost_bound}.

    Finally, we note that $|C^c| = (1-\beta)m$ and $f_t(|\langle \va_{j}, \ve^{(t-1)}\rangle|) \ge f_t(\|\ve^{(t-1)}\|)$, which follows again from the nonincreasing property of $f_t$ and~\eqref{eq:CS_unitnorm}.  Thus, we have $$\sum_{j\not\in C} f_t(|\langle \va_{j}, \ve^{(t-1)}\rangle|) \ge (1-\beta)m f_t(\|\ve^{(t-1)}\|),$$
    which upper bounds the last term, obtaining~\eqref{eq:main-bound}.
\end{proof}

In order to better explore Theorem~\ref{thm:nonincreasing_f_main}, we specialize it to three types of weighting functions in Section~\ref{subsec:special_cases}: (1) threshold weighting functions, which uniformly sample from the rows with residual magnitudes below a threshold value; (2) exponential weighting functions, which employ a function that exponentially decreases with the residual magnitudes; and (3) inverse weighting functions, which decrease proportionally to the inverse of a power of the residual magnitude.  We explore these function types in Subsections~\ref{subsubsec:threshold_type},~\ref{subsubsec:exponential_type}, and~\ref{subsubsec:inverse_type}, respectively.  In each section, we specialize Theorem~\ref{thm:nonincreasing_f_main} to a corollary providing a single-iteration bound on the expected error of Algorithm~\ref{alg:gwrk}, GWRK, when using a weighting function $f_t$ of this type.

\subsection{Special Cases}\label{subsec:special_cases}

\begin{table}[htbp]
    \centering

    \renewcommand{\tabularxcolumn}[1]{m{#1}}

    \setlength{\weightingTableWidth}{
        \dimexpr\linewidth-6\tabcolsep\relax
    }

    \newcommand{\wtTypeWidth}{0.16\weightingTableWidth}
    \newcommand{\wtFormWidth}{0.31\weightingTableWidth}
    \newcommand{\wtExampleWidth}{0.32\weightingTableWidth}
    \newcommand{\wtDescWidth}{0.19\weightingTableWidth}

    \begin{tabularx}{\linewidth}{
        @{}
        >{\raggedright\arraybackslash}m{\wtTypeWidth}
        >{\raggedright\arraybackslash}m{\wtFormWidth}
        >{\raggedright\arraybackslash}X
        >{\raggedright\arraybackslash}m{\wtDescWidth}
        @{}
    }
        \toprule
        \textbf{Weighting Type}
        & \textbf{Form}
        & \textbf{Example Function}
        & \textbf{Description}
        \\
        \midrule

        \textbf{Threshold Weighting}\newline
        Section~\ref{subsubsec:threshold_type}
        &
        $\displaystyle
        f_t(r)=
        \begin{cases}
            1, & r\le\theta^{(t-1)}, \\[2pt]
            0, & r>\theta^{(t-1)}.
        \end{cases}$
        &
        \multicolumn{2}{l@{}}{            \begin{tabular}[c]{
                @{}
                >{\raggedright\arraybackslash}m{\wtExampleWidth}
                >{\raggedright\arraybackslash}m{\wtDescWidth}
                @{}
            }
                $\theta^{(t-1)}\!=\!Q_{\tau}(|\vr^{(t-1)}|)$
                &
                Quantile threshold (QRK)
                \\
                \cmidrule(l){1-2}

                $\theta^{(t-1)}
                \!=\!\theta\!=\!Q_{\scriptscriptstyle 0.5}(|\vr^{(0)}|)$
                &
                Constant threshold
            \end{tabular}        }
        \\
        \midrule

        \textbf{Exponential\newline Weighting}\newline
        Section~\ref{subsubsec:exponential_type}
        &
        $\displaystyle
        f_t(r)\!=\!
        \exp\!\left(
            -\frac{r^2}{g(|\vr^{(t-1)}|)}
        \right)$
        &
        \multicolumn{2}{l@{}}{            \begin{tabular}[c]{
                @{}
                >{\raggedright\arraybackslash}m{\wtExampleWidth}
                >{\raggedright\arraybackslash}m{\wtDescWidth}
                @{}
            }
                $g(|\vr^{(t-1)}|)
                \!=\!Q_{\scriptscriptstyle 0.5}(|\vr^{(t-1)}|)$
                &
                Exponential median
                \\
                \cmidrule(l){1-2}

                $g(|\vr^{(t-1)}|)
                \!=\!Q^2_{\scriptscriptstyle 0.5}(|\vr^{(t-1)}|)$
                &
                Exponential squared median
                \\
                \cmidrule(l){1-2}

                $g(|\vr^{(t-1)}|)
                \!=\!\sigma_{\scriptscriptstyle\mathrm{MAD}}
                (|\vr^{(t-1)}|)$
                &
                Exponential $\sigma_{\scriptscriptstyle\mathrm{MAD}}$
                \\
                \cmidrule(l){1-2}

                $g(|\vr^{(t-1)}|)
                \!=\!\sigma_{\scriptscriptstyle\mathrm{MAD}}^2
                (|\vr^{(t-1)}|)$
                &
                Exponential $\sigma_{\scriptscriptstyle\mathrm{MAD}}^2$
            \end{tabular}        }
        \\
        \midrule

        \textbf{Inverse Weighting}\newline
        Section~\ref{subsubsec:inverse_type}
        &
        $\displaystyle
        f_t(r)=
        \begin{cases}
            r^{-q}, & r>0, \\[2pt]
            0, & r\le 0.
        \end{cases}$
        &
        \multicolumn{2}{l@{}}{            \begin{tabular}[c]{
                @{}
                >{\raggedright\arraybackslash}m{\wtExampleWidth}
                >{\raggedright\arraybackslash}m{\wtDescWidth}
                @{}
            }
                $q=1.0$
                &
                Inverse absolute magnitude
                \\
                \cmidrule(l){1-2}

                $q=1.2$
                &
                --
                \\
                \cmidrule(l){1-2}

                $q=2.0$
                &
                Inverse squared magnitude
            \end{tabular}        }
        \\
        \bottomrule
    \end{tabularx}
    \caption{Three special weighting families and the corresponding examples.}
    \label{tab:mytable}
\end{table}

\begin{figure}[htbp]
    \centering
    \includegraphics[width=\linewidth]{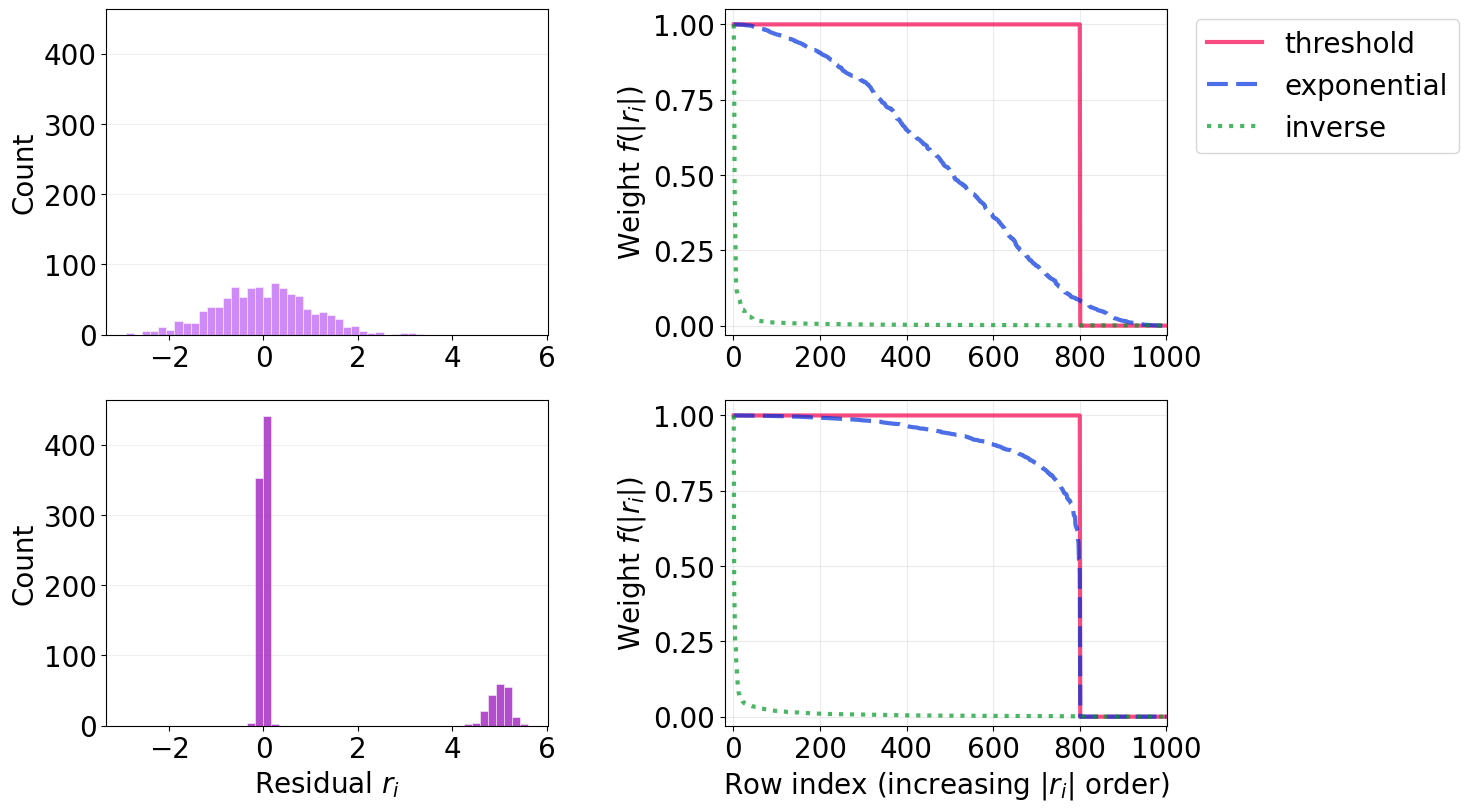}
    \caption{
    The top row shows residuals without a clear separation in magnitude (single peak), while
    the bottom row shows well-separated residuals (two peaks).
    The left panels display histograms of the residuals $r_i$;
    the right panels display the corresponding weights, with entries ordered
    by increasing absolute residual $|r_i|$.
    The weighting function families are:
    threshold (solid red), exponential (dashed blue), and inverse (dotted green).
    Inverse weights are normalized by their maximum for visualization.
    }
    \label{fig:wfshape}
\end{figure}

We now examine more closely the behavior of GWRK, Algorithm~\ref{alg:gwrk}, and the convergence guarantee, Theorem~\ref{thm:nonincreasing_f_main}, when the weighting function $f_t$ is one of three special function types.  We summarize these function types and give several examples of each in Table~\ref{tab:mytable}, and visualize an archetypal example of each type of weighting function in Figure~\ref{fig:wfshape}.  We next discuss these three function types, and some advantages and disadvantages of each.  

The first special weighting function type is \emph{threshold} weighting functions, which uniformly weight the rows corresponding to residual entries with magnitude below a given (fixed or adaptive) threshold value; see Subsection~\ref{subsubsec:threshold_type} for more details.  This function type encompasses Quantile Randomized Kaczmarz (QRK) where the threshold is determined by the empirical quantile of the residual magnitudes.  We note that this type of function likely requires some estimate of the magnitude of corruptions, or, like in the case of QRK, a good estimate of the fraction of equations that are corrupted, $\beta$, in order to produce a threshold that correctly filters out corrupted residual entries as intended.  On the other hand, though, one is able to provide a thorough convergence analysis -- indeed,~\cite{haddock2022quantile,steinerberger2023quantile} prove that QRK converges to the solution of the uncorrupted system linearly in expectation under mild assumptions. 

The second special weighting function type is \emph{exponential} weighting functions, which exponentially downweight rows according to their corresponding residual magnitude.  Unlike the threshold weighting functions, and especially quantile-based functions, these methods do not require an estimate of the fraction of corrupted equations -- in this way, we say that these methods are \emph{parameter-free}.  However, without this assumed knowledge about the corruption fraction $\beta$ or corruption magnitudes, these functions must necessarily be more conservative than the uniform thresholded functions and assign more probability mass to rows corresponding to smaller magnitude residual entries; see Figure~\ref{fig:wfshape} for a visualization of the conservative nature of these functions.  

The last and most conservative special weighting function type is \emph{inverse} weighting functions, which assign rows weights inversely proportional to the magnitude of their corresponding residual entries.  These methods are also parameter-free, and require no knowledge of the corruption fraction $\beta$ or corruption magnitudes.  However, these methods are generally the most conservative type of function considered -- they assign nearly all of the probability mass to very few rows corresponding to the smallest residual entries; see Figure~\ref{fig:wfshape}.  

Figure~\ref{fig:wfshape} demonstrates the shape of weighting function $f_t$
given two different distributions of residuals $r_i$. In the top row, the residuals are normally distributed (see top-left subplot), and residuals corresponding to ``corruptions" are not obvious based on residual values. In this setting, we see that the threshold weighting function (pink solid line) assigns zero weight to the 20\% largest in magnitude residuals, the exponential function's weights slowly decay (blue dashed line), and the inverse function's weights concentrate on the indices corresponding to the smallest residuals (green dotted line). The bottom row repeats this experiment except with separation in the residuals (large residuals correspond exactly to the corrupted rows). In this case, the inverse-type weighting still concentrates weights on small residuals, whereas the exponential-type weights any residual that is not too large, beginning to approximate the weights assigned by the threshold-type. 
Lastly, the threshold weighting function remains the same under the two different distributions of residuals.

\subsubsection{Thresholded Weighting Functions}\label{subsubsec:threshold_type}

We now turn our focus to the thresholded weighting functions.  With an assumed bound on the threshold in terms of the error, we prove a single-iteration bound on the expected error of GWRK with this type of weighting function in Corollary~\ref{cor:thresholded_convergence}.  Under mild assumptions, we can iterate this bound to prove that the expected squared error converges to zero linearly, with no convergence horizon. 

Throughout what follows, we assume that $\tau_t m$ and $\beta m$ are integers. Given an integer $l$, define 
\begin{equation}
    s_{l,\min}^2(\mA) = \min_{S \subset [m], |S| = l} s_{\min}^2(\mA_S).
    \label{eq:sigma_betamin}
\end{equation}

\begin{corollary}[Thresholded convergence]
Suppose $\vb = \mA\vx^* + \vc$ where $\mA \in \mathbb{R}^{m \times n}$ is row-normalized and $C = \text{supp}(\vc)$ satisfies $|C| = \beta m$.
Let $\vx^{(t)}$ denote the iterates produced in the $t$th iteration of Algorithm~\ref{alg:gwrk}, where  
$$f_t(r) = \begin{cases}
            1 & \text{ if } r \le \theta^{(t-1)} \\ 0 & \text{ if } r > \theta^{(t-1)}
        \end{cases}.$$ Suppose that $\left|\left\{j \in [m] : f_t(|r_j^{(t-1)}|) = 1\right\}\right| \ge \tau_t m$ with $1 > \tau_t > \beta$ and that there exists $\kappa_t > 0$ such that  
        $$\theta^{(t-1)} \le \kappa_t \|\vx^{(t-1)} - \vx^*\|.$$ 
        Define $N_t = (\tau_t - \beta)m$.  Then the iterates satisfy
\[
  \mathbb{E}_{t-1}\big[\|\vx^{(t)}-\vx^*\|^2\big]
  \;\le\;
  \left(
    1 \;-\; \frac{s^2_{N_t,\min}(\mA)}{m}
    \;+\; \frac{\beta}{\tau_t-\beta}\big(\kappa_t^2 + 2\kappa_t\big)
  \right) \|\vx^{(t-1)}-\vx^*\|^2 ,
\]
where $s^2_{N_t,\min}(\mA)$ is as defined in~\eqref{eq:sigma_betamin}. In particular, if
\[
  1 - \frac{s^2_{N_t,\min}(\mA)}{m} + \frac{\beta}{\tau_t-\beta}\big(\kappa_t^2+2\kappa_t\big) < 1 - \epsilon
\]
for some fixed $0 < \epsilon < 1$ for all $t$, then
$\mathbb{E}\|\vx^{(t)}-\vx^*\|^2 \to 0$ geometrically, with no convergence
horizon. \label{cor:thresholded_convergence}
\end{corollary}

\begin{remark}
    The Quantile Randomized Kaczmarz method, which is GWRK with $$f_t(r) = \begin{cases}
            1 & \text{ if } r \le \theta^{(t-1)} \\ 0 & \text{ if } r > \theta^{(t-1)}
        \end{cases}$$ where $\theta^{(t-1)} = Q_\tau(|\vr^{(t-1)}|)$ and $\beta < \tau < 1-\beta$ with $\beta < 0.5$, obeys Corollary~\ref{cor:thresholded_convergence} with constants $\tau_t = \tau$, $N_t = N := (\tau-\beta)m$, and 
        \begin{equation}
        \kappa_t = \kappa := \frac{s_{\max}(\mA)}{\sqrt{m}\sqrt{1-\tau-\beta}}.
        \label{eq:kappa}
        \end{equation}
        This value for $\kappa$ is given in~\cite[Lemma 1]{steinerberger2023quantile}.  More generally, Corollary~\ref{cor:thresholded_convergence} applies to a variant of QRK in which the quantile level $\tau_t$ changes across iterations or is adaptive, as long as $1-\beta > \tau_t > \beta$ is maintained.  In this case, Corollary~\ref{cor:thresholded_convergence} holds with $\kappa_t := \frac{s_{\max}(\mA)}{\sqrt{m}\sqrt{1-\tau_t-\beta}}$. 
\end{remark}

\begin{remark}
    We note that Corollary~\ref{cor:thresholded_convergence} applied to GWRK with a fixed threshold weighting function as described in Table~\ref{tab:mytable} is not able to provide a convergence horizon-free result.  In particular, the assumption $\theta^{(t-1)} = \theta \le \kappa_t \|\ve^{(t-1)}\|$ would require increasing $\kappa_t \rightarrow \infty$ as $\|\ve^{(t-1)}\| \rightarrow 0$.  However, the single iteration bound holds with $\kappa_t$ chosen such that $\kappa_t > \frac{\theta}{\|\ve^{(t-1)}\|}.$
\end{remark}

We now provide a proof for Corollary~\ref{cor:thresholded_convergence}.
\begin{proof} 
    Note that $f_t(0) = 1$. We first lower-bound the denominator
$\sum_{j\notin C} f_t(|\langle \va_j, \vx^{(t-1)}-\vx^*\rangle|)$ appearing in~\eqref{eq:almost_bound} 
by a counting argument. For $j \notin C$ we have $\vc_j = 0$, so
$\langle \va_j, \vx^{(t-1)}-\vx^*\rangle = \vr_j^{(t-1)}$, and thus
\[
  \sum_{j\notin C} f_t\big(|\langle \va_j, \vx^{(t-1)}-\vx^*\rangle|\big)
  \;=\; \big|\{\,j \notin C : f_t(|r_j^{(t-1)}|) = 1\}\big| .
\]
By the assumption that there are at least $\tau_t m$ indices of $j$ such that $f_t(|r_j^{(t-1)}|) = 1$, and since at most $|C| = \beta m$
of these can lie in $C$, we have that at least $(\tau_t-\beta)m$ of them
lie in $C^c$. Hence
\begin{equation}
  \sum_{j\notin C} f_t\big(|\langle \va_j, \vx^{(t-1)}-\vx^*\rangle|\big)
  \;\ge\; (\tau_t-\beta)m =N_t ,\label{eq:prob_sum_bound_threshold}
\end{equation}
which is strictly positive precisely because $\tau_t > \beta$. 

Now, we turn to bounding $s_{\min}^2(\tilde{\mA}_{C^c})$ in~\eqref{eq:almost_bound}.  Noting that $f_t$ takes only values in $\{0,1\}$, we note that $\tilde{\mA}_{C^c}$ has only rows identical to those of $\mA_{C^c}$ or rows of all zeros.  Furthermore, the rows that remain the same are those for which the corresponding value of $f_t(|r_j^{(t-1)}|) = 1$.    By the argument before, we know that the size of this set is at least $N_t$.  Thus, we have $$s_{\min}^2(\tilde{\mA}_{C^c}) \ge s^2_{N_t,\min}(\mA).$$

Substituting this,
$f_t(0)=1$, and the bound~\eqref{eq:prob_sum_bound_threshold} into~\eqref{eq:almost_bound}, we have
\[
  \mathbb{E}_{t-1}\big[\|\vx^{(t)}-\vx^*\|^2\big]
  \le
  \left(1 - \frac{s^2_{N_t,\min}(\mA)}{m}\right)\|\vx^{(t-1)}-\vx^*\|^2
  + \frac{1}{N_t}
  \sum_{j\in C_F^{(t)}} \big[c_j^2 - |\langle \va_j, \vx^{(t-1)}-\vx^*\rangle|^2\big] .
\]

Now, we note that $$|c_j| = |\langle \va_{j}, \vx^{(t-1)} - \vx^*\rangle - c_j - \langle \va_{j}, \vx^{(t-1)} - \vx^*\rangle| \le |r_j^{(t-1)}| + |\langle \va_{j}, \vx^{(t-1)} - \vx^*\rangle|$$ and thus, by squaring both sides, we have $$c_j^2 - |\langle \va_{j}, \vx^{(t-1)} - \vx^*\rangle|^2 \le |r_j^{(t-1)}|^2 + 2 |r_j^{(t-1)}||\langle \va_{j}, \vx^{(t-1)} - \vx^*\rangle| \le |r_j^{(t-1)}|^2 + 2 |r_j^{(t-1)}| \|\vx^{(t-1)} - \vx^*\|,$$ 
where the last inequality follows from~\eqref{eq:CS_unitnorm}.
By definition, $j \in C_F^{(t)}$ requires $f_t(|r_j^{(t-1)}|) > 0$, i.e.\
$|r_j^{(t-1)}| \le \theta^{(t-1)}$. Since $|C_F^{(t)}| \le |C| = \beta m$,
\[
  \sum_{j\in C_F^{(t)}} \big[c_j^2 - |\langle \va_j,\vx^{(t-1)}-\vx^*\rangle|^2\big]
  \;\le\; \beta m\big((\theta^{(t-1)})^2 + 2\theta^{(t-1)}\|\vx^{(t-1)}-\vx^*\|\big) .
\]
By assumption,
$\theta^{(t-1)} \le \kappa_t\,\|\vx^{(t-1)}-\vx^*\|$, so
\[
  \sum_{j\in C_F^{(t)}} \big[c_j^2 - |\langle \va_j,\vx^{(t-1)}-\vx^*\rangle|^2\big]
  \;\le\; \beta m\big(\kappa_t^2 + 2\kappa_t\big)\|\vx^{(t-1)}-\vx^*\|^2 .
\]
Thus, putting it all together, we have
\[
  \mathbb{E}_{t-1}\big[\|\vx^{(t)}-\vx^*\|^2\big]
  \le
  \left(1 - \frac{s^2_{N_t,\min}(\mA)}{m}
  + \frac{\beta}{\tau_t - \beta}\big(\kappa_t^2+2\kappa_t
  \big)\right)
  \|\vx^{(t-1)}-\vx^*\|^2,
\]
which is the claimed bound.  Assuming that $\left(1 - \frac{s^2_{N_t,\min}(\mA)}{m}
  + \frac{\beta}{\tau_t - \beta}\big(\kappa_t^2+2\kappa_t
  \big)\right) < 1 -\epsilon$ for all $t$, we can iterate with the tower of expectation to guarantee $\mathbb{E}[\|\ve^{(t)}\|^2] \rightarrow 0$.
\end{proof}

\subsubsection{Exponential Weighting Functions}\label{subsubsec:exponential_type}

\begin{figure}
    \centering
    \includegraphics[width=0.48\linewidth]{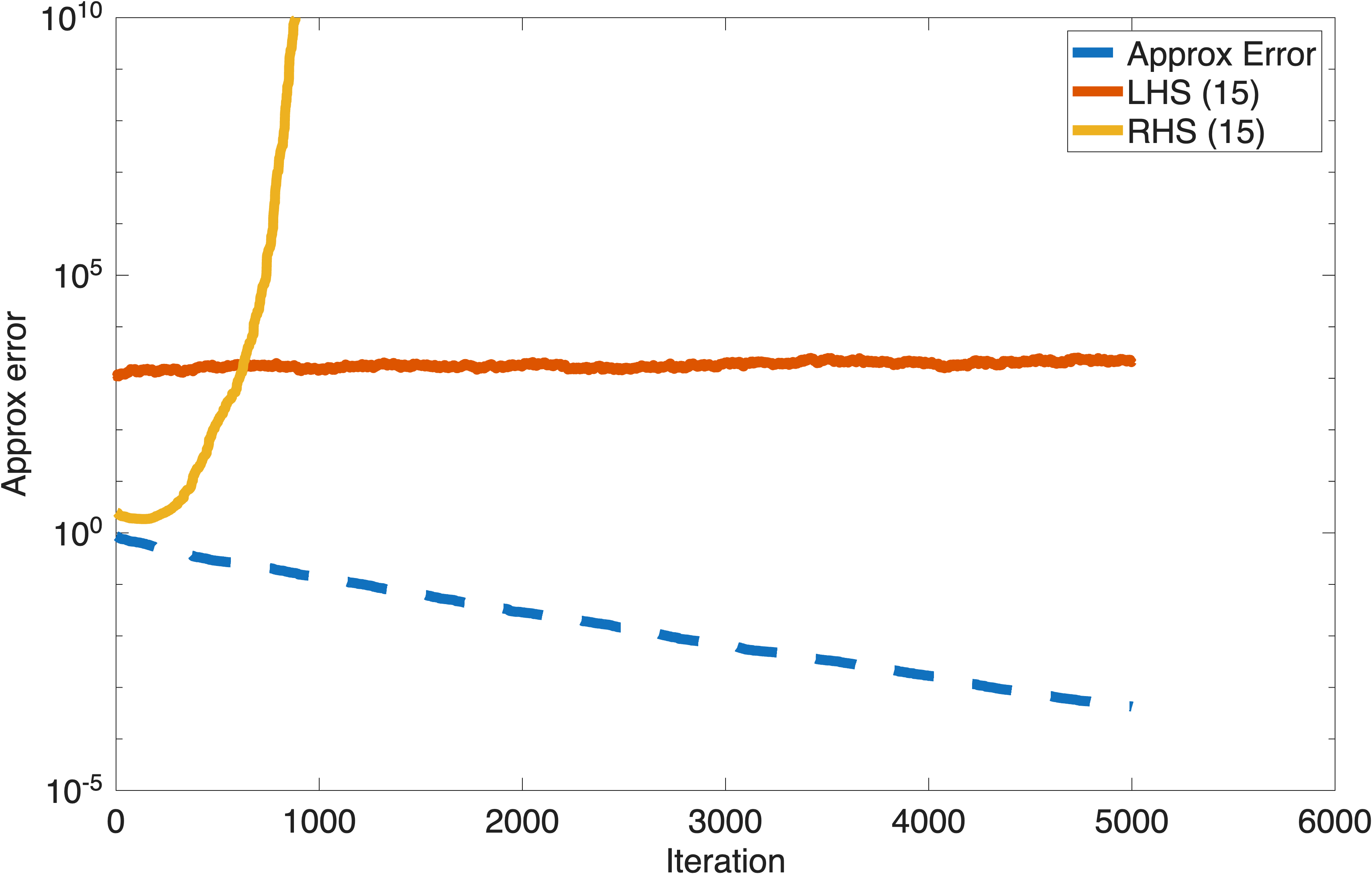}
    \includegraphics[width=0.48\linewidth]{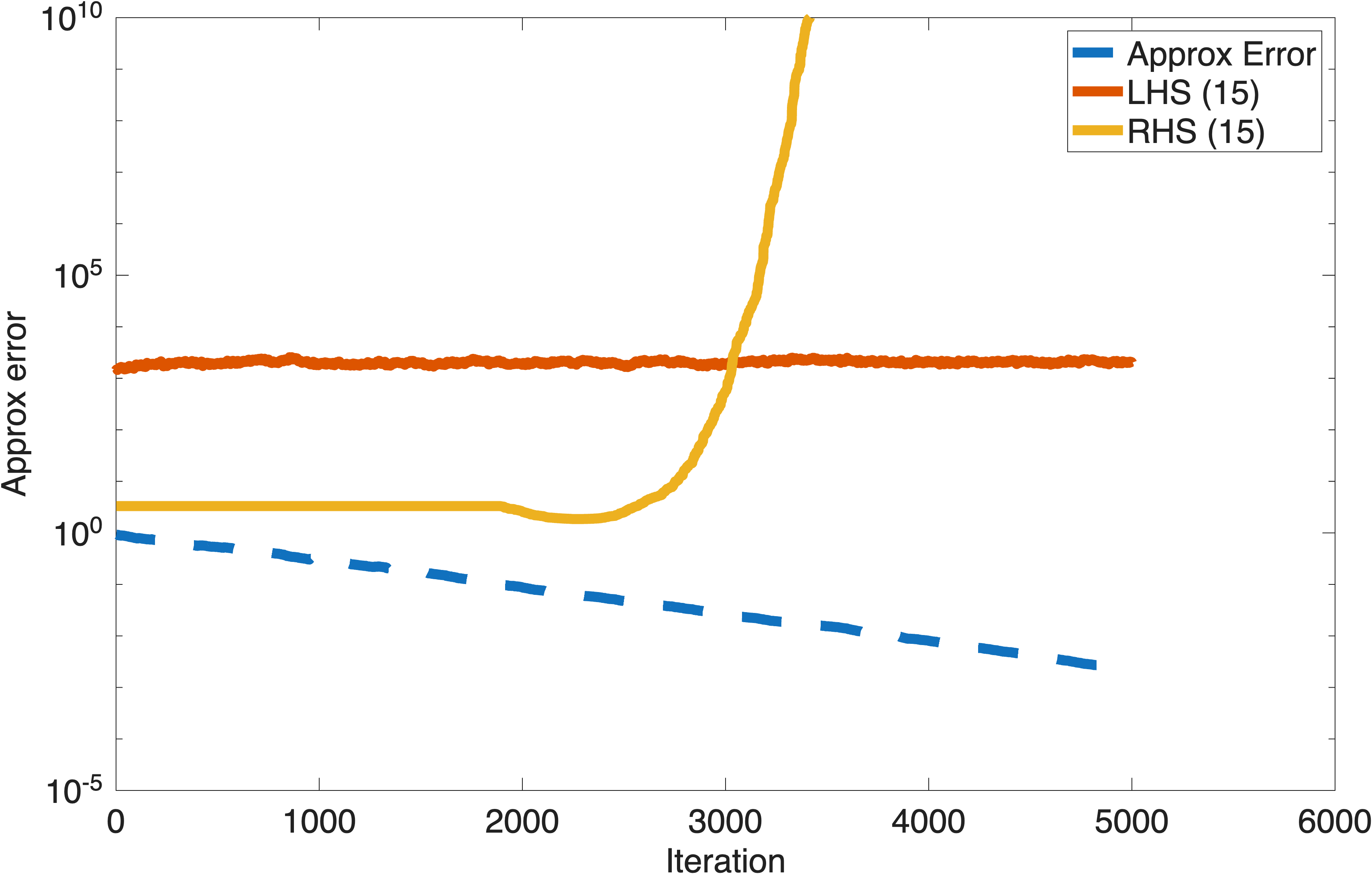} \\
    \includegraphics[width=0.48\linewidth]{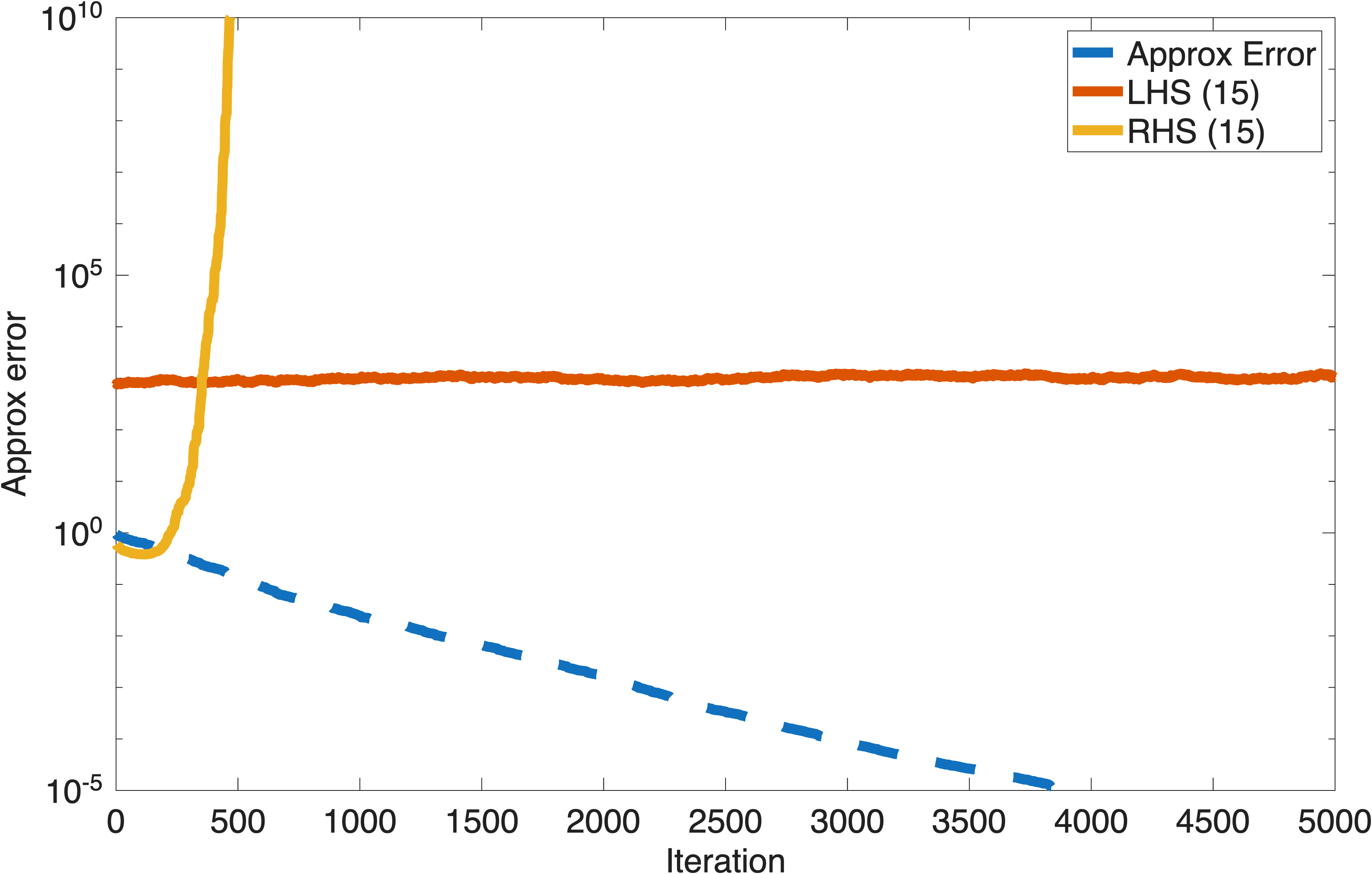}
    \includegraphics[width=0.48\linewidth]{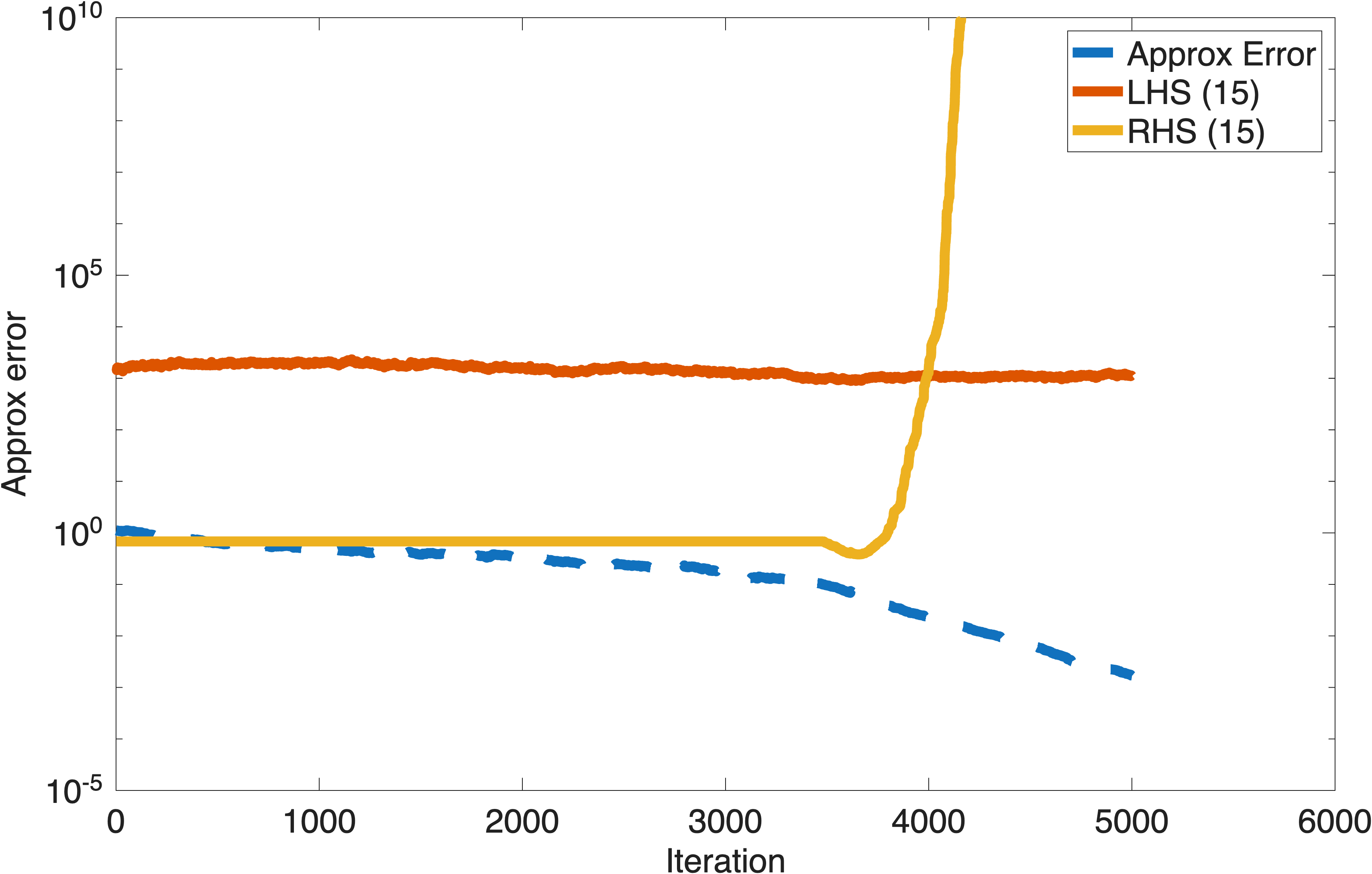} \\
    \includegraphics[width=0.48\linewidth]{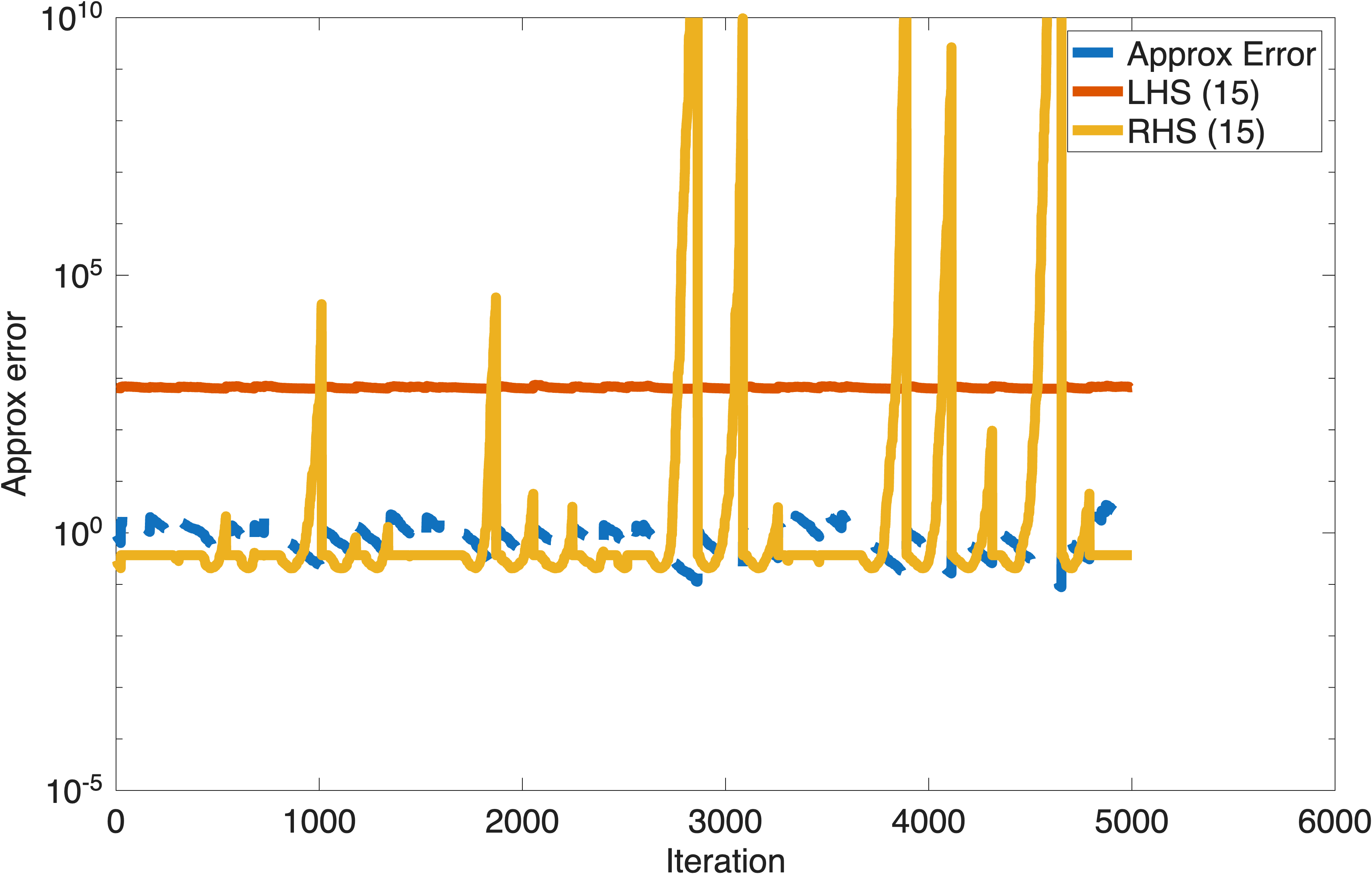}
    \includegraphics[width=0.48\linewidth]{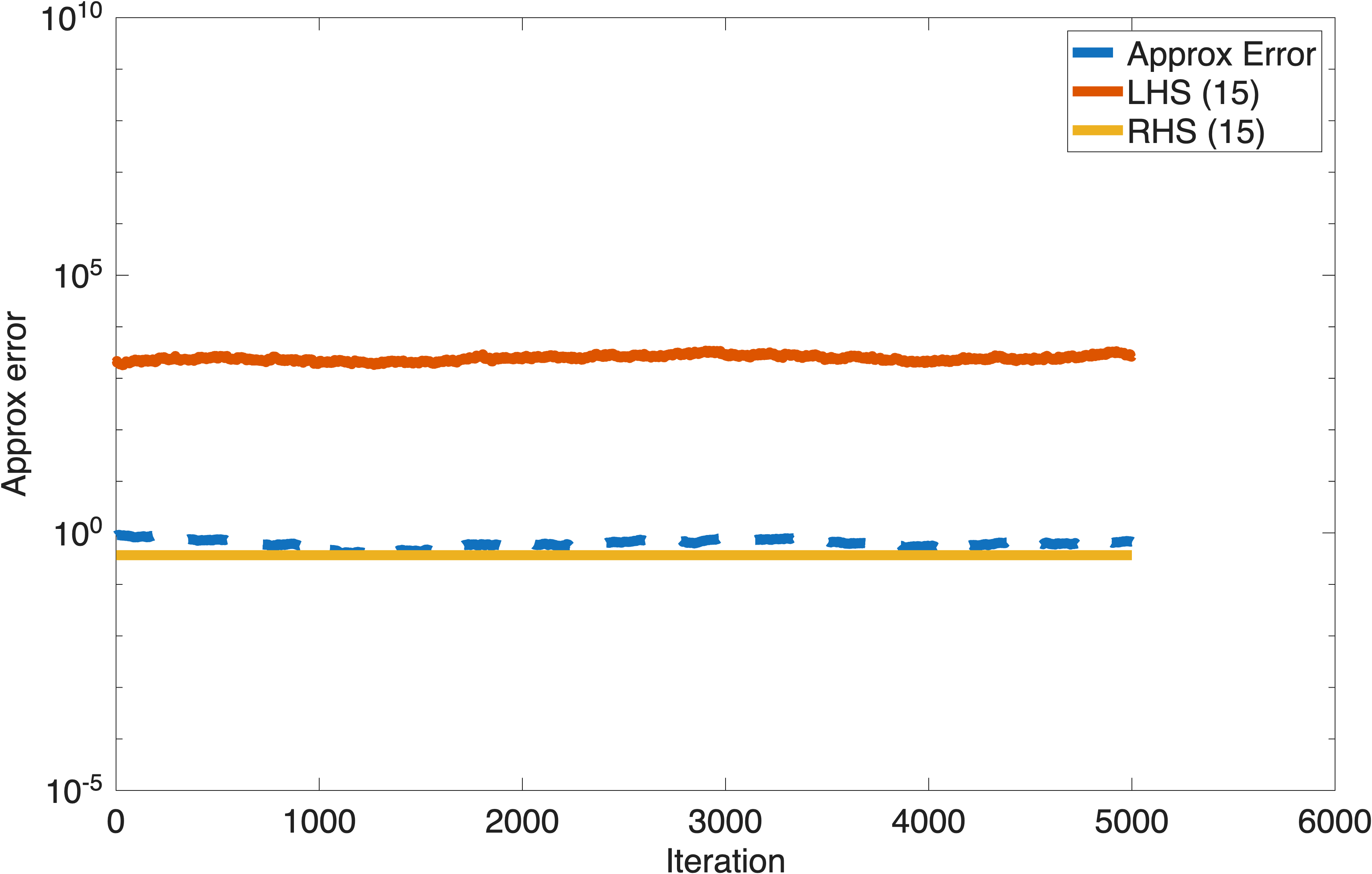}
    \caption{$\beta = 0.10$ (top), $\beta = 0.35$ (middle), $\beta = 0.5$ (bottom), with the corruption values $c_j = 1$ (left column) and $c_j = 0.1$ (right column). Condition~\eqref{eq:medexp_condition} requires $LHS < RHS$.}
    \label{fig:medexp_condition}
\end{figure}

We next turn our focus to the exponential weighting functions.  With an assumed bound on the parametrizing function $g$ in terms of the error, we prove a single iteration bound on the expected error of GWRK with this type of weighting function in Corollary~\ref{cor:medianexp}.  Unlike Corollary~\ref{cor:thresholded_convergence} with the thresholded weighting functions, it is not clear that one can iterate this bound and form a convergence guarantee with no convergence horizon.  We explore this possibility under a simplifying assumption in Remark~\ref{rem:median_conv_horizon} and Figure~\ref{fig:medexp_condition}.  

\begin{corollary}[Exponential convergence]
Suppose $\vb = \mA\vx^* + \vc$ where $\mA \in \mathbb{R}^{m \times n}$ is row-normalized and $C = \text{supp}(\vc)$ satisfies $|C| = \beta m$ with $\beta < 0.5$.
Let $\vx^{(t)}$ denote the iterates produced in the $t$th iteration of Algorithm~\ref{alg:gwrk}, where
\begin{equation}
    f_t(r) = \exp(-r^2/g(|\vr^{(t-1)}|)),
    \label{eq:ff_medexp}
\end{equation} and $g$ is a function satisfying $0 < g(|\vr^{(t-1)}|) \le \kappa \|\vx^{(t-1)} - \vx^*\|^\alpha$ with $\kappa > 0$ and $\alpha > 0$.  
We define $\ve^{(t)} = \vx^{(t)} - \vx^*$ to be the error in the $t^{th}$ iteration of Algorithm~\ref{alg:gwrk}.  Define $c_{\min}^{(t-1)} = \max\{\min_{j \in C} |c_j|, \|\ve^{(t-1)}\|\}$ and let $\tilde{\mA}_{C^c}$ be the matrix with $j$th row $\tilde{\va}_j := \sqrt{f_t(|r_j^{(t-1)}|)} \va_j$. Then if $c_{\min}^{(t-1)} - 2\|\ve^{(t-1)}\| > 0$, the iterates of Algorithm~\ref{alg:gwrk} satisfy
\begin{align*}
  \mathbb{E}_{t-1}\big[\|\ve^{(t)}\|^2\big]
  \;&\le\;  \left(1 - \frac{s_{\min}^2(\tilde{\mA}_{C^c})}{m}\right) \|\ve^{(t-1)}\|^2 + \frac{\|\vc\|^2}{(1-\beta)m} \exp \left(-\frac{c_{\min}^{(t-1)}(c_{\min}^{(t-1)} - 2\|\ve^{(t-1)}\|)}{\kappa\|\ve^{(t-1)}\|^\alpha}\right) \\
  &\le\;  \left(1 - \frac{s_{\min}^2(\tilde{\mA}_{C^c})}{m}\right) \|\ve^{(t-1)}\|^2 + \frac{2\|\vc\|^2}{m} \exp \left(-\frac{c_{\min}^{(t-1)}(c_{\min}^{(t-1)} - 2\|\ve^{(t-1)}\|)}{\kappa\|\ve^{(t-1)}\|^\alpha}\right).
\end{align*}
\label{cor:medianexp}
\end{corollary}

\begin{remark}
We note that we consider several particular functions $g$, and for each of these, we have a value $\kappa$ and $\alpha$ that suffice for the assumed bound $g(|\vr^{(t-1)}|) \le \kappa \|\vx^{(t-1)} - \vx^*\|^\alpha$.  These are stated and proven in Lemma~\ref{lem:exponential_function_bounds}.
\end{remark}

We now prove Corollary~\ref{cor:medianexp}.

\begin{proof} 
    Beginning from~\eqref{eq:almost_bound}, we note that $f_t(0) = 1$,
 and bound the second term as in~\eqref{eq:main-bound}, which gives   \begin{align*}
        \mathbb{E}_{t-1}[\|\ve^{(t)}\|^2] &\le \left(1 - \frac{s_{\min}^2(\tilde{\mA}_{C^c})}{m}\right) \|\ve^{(t-1)}\|^2 +\frac{f_t(c_{\min}^{(t-1)} - \|\ve^{(t-1)}\|)}{(1 - \beta)mf_t(\|\ve^{(t-1)}\|)} \|\vc\|^2\\
        &\le \left(1 - \frac{s_{\min}^2(\tilde{\mA}_{C^c})}{m}\right) \|\ve^{(t-1)}\|^2 + \frac{\|\vc\|^2}{(1 - \beta)m} \exp \left(-\frac{c_{\min}^{(t-1)}(c_{\min}^{(t-1)} - 2\|\ve^{(t-1)}\|)}{\kappa\|\ve^{(t-1)}\|^\alpha}\right) 
    \end{align*}
  where $f_t(c_{\min}^{(t-1)} - \|\ve^{(t-1)}\|)/f_t(\|\ve^{(t-1)}\|)$ in the last term   has been simplified using exponential properties. Thus, the bound follows from the assumed bounds $g(|\vr^{(t-1)}|) \le \kappa\|\ve^{(t-1)}\|^\alpha$ and $\beta < 0.5.$

\end{proof}

Now, in order to prove bounds of the form $g(|\vr^{(t-1)}|) \le \kappa \|\vx^{(t-1)} - \vx^*\|^\alpha$ for example functions $g$, we need an additional lemma, which yields bounds on the median and $\sigma_{\text{MAD}}$ for input vectors with a majority of entries of bounded magnitude.

\begin{lemma}
    Suppose vector $\vv$ has more than half of entries in $[-d,d]$.  Then $$Q_{0.5}(|\vv|) \le d$$ and $$Q_{0.5}(||\vv| - Q_{0.5}(|\vv|)|) \le d.$$ \label{lem:median_bounds}
\end{lemma}
\begin{proof}
    Since more than half of the entries of $\vv$ have magnitude at most $d$, there exists some index $j$ where $|v_j| \le d$ and $Q_{0.5}(|\vv|) \le |v_j|$, so the first bound holds.  

    For the second bound, note that by the same argument, there exists some $j$ where $|v_j| \le d$ and $Q_{0.5}(||\vv| - Q_{0.5}(|\vv|)|) \le ||v_j| - Q_{0.5}(|\vv|)|$.  Now, both $|v_j|$ and $Q_{0.5}(|\vv|)$ lie in $[0,d]$, so their difference has magnitude no larger than $d$.  
\end{proof}

With this lemma in hand, we can now provide bounds of the form $g(|\vr^{(t-1)}|) \le \kappa \|\vx^{(t-1)} - \vx^*\|^\alpha$ for all of the $g$ functions utilized in the example functions of exponential weighting type in Table~\ref{tab:mytable}.  

\begin{lemma}
Assume the setting of Corollary~\ref{cor:medianexp} holds.  Then we have the following bounds:
\begin{enumerate}
    \item If $g(|\vr^{(t-1)}|) = Q_{0.5}(|\vr^{(t-1)}|)$, then $g(|\vr^{(t-1)}|) \le \|\vx^{(t-1)} - \vx^*\|.$
    \item If $g(|\vr^{(t-1)}|) = Q^2_{0.5}(|\vr^{(t-1)}|)$, then $g(|\vr^{(t-1)}|) \le \|\vx^{(t-1)} - \vx^*\|^2.$
    \item If $g(|\vr^{(t-1)}|) = \sigma_{\text{MAD}}(|\vr^{(t-1)}|)$, then $g(|\vr^{(t-1)}|) \le \frac32 \|\vx^{(t-1)} - \vx^*\|.$
    \item If $g(|\vr^{(t-1)}|) = \sigma_{\text{MAD}}^2(|\vr^{(t-1)}|)$, then $g(|\vr^{(t-1)}|) \le \frac94 \|\vx^{(t-1)} - \vx^*\|^2.$
\end{enumerate}\label{lem:exponential_function_bounds}
\end{lemma}

\begin{proof}
    Now, note that since $\beta < 0.5$, more than half of the entries in $\vr^{(t-1)}$ correspond to uncorrupted entries.  For any $j \not\in C$, we have $$|r_j^{(t-1)}| = |\langle \va_{j}, \vx^{(t-1)} - \vx^*\rangle| \le \|\vx^{(t-1)} - \vx^*\|$$ and thus, by Lemma~\ref{lem:median_bounds} we have $Q_{0.5}(|\vr^{(t-1)}|) \le \|\vx^{(t-1)} - \vx^*\|$ and \[\sigma_{\text{MAD}}(|\vr^{(t-1)}|) = \frac{1}{0.6745} Q_{0.5}(||\vr^{(t-1)}| - Q_{0.5}(|\vr^{(t-1)}|)|) \le \frac{1}{0.6745} \|\vx^{(t-1)} - \vx^*\|.\]  All four bounds follow.
\end{proof}

\begin{remark}
    We note that Corollary~\ref{cor:medianexp} holds only when $c_{\min}^{(t-1)} > 2\|\ve^{(t-1)}\|$, and thus, only once the residual values have become sorted with corrupted residual magnitudes strictly larger than uncorrupted residual magnitudes.  In order to relax this assumption, one would need a lower bound on $g(|\vr^{(t-1)}|)$ of the form $\kappa_l \|\ve^{(t-1)}\|^\alpha$.  Such a bound is not possible to provide for $\sigma_{\text{MAD}}$ or $\sigma_{\text{MAD}}^2$ in general.  However, one can prove that $$Q_{0.5}^2(|\vr^{(t-1)}|) \ge \frac{2}{m} s_{l,\min}^2(\mA) \|\ve^{(t-1)}\|^2,$$ where $l = (\tfrac12 - \beta)m.$  Thus, a bound even in this early convergence regime is possible for the exponential weighting functions where $g$ is defined by the (squared or unsquared) median of the residual.  However, we do not pursue this bound in detail since it does not provide a convergence-horizon free guarantee in general. \label{rem:med_lower_bound}
\end{remark}

\begin{remark}\label{rem:median_conv_horizon}
    A natural question is whether Corollary~\ref{cor:medianexp} can offer a bound without convergence horizon.  To investigate, we make a simplifying assumption and reparameterize the error.  Suppose that the corruption vector has all nonzero coordinates with equal magnitude, so $\|\vc\|^2 \le \beta m (c_{\min}^{(t-1)})^2$.  We also introduce the parameter $\gamma_t$ and let $\|\vx^{(t-1)} - \vx^*\| = \gamma_t c_{\min}^{(t-1)}$. In this case, we have \[
  \mathbb{E}_{t-1}\big[\|\vx^{(t)}-\vx^*\|^2\big]
  \;\le\;  \left( \gamma_t^2 - \frac{\gamma_t^2s_{\min}^2(\tilde{\mA}_{C^c})}{m} +  \frac{\beta}{1-\beta}\exp \left(-\frac{1-2\gamma_t}{\kappa\gamma_t^\alpha (c_{\min}^{(t-1)})^{\alpha - 2}} \right) \right) \Big(c_{\min}^{(t-1)}\Big)^2.
\] 
In order to guarantee that $\mathbb{E}_{t-1}\big[\|\vx^{(t)}-\vx^*\|^2\big] < \|\vx^{(t-1)} - \vx^*\|^2 = \gamma_t^2 (c_{\min}^{(t-1)})^2$, it suffices to have  

\begin{equation}
     \frac{m}{s_{\min}^2(\tilde{\mA}_{C^c})} < \frac{\gamma_t^2 (1-\beta)}{\beta} \exp \left(\frac{1 - 2\gamma_t}{\kappa\gamma_t^\alpha (c_{\min}^{(t-1)})^{\alpha - 2}}\right).
    \label{eq:medexp_condition}
\end{equation}
\end{remark}

We now test the viability of condition~\eqref{eq:medexp_condition}. To do so, we run Algorithm~\ref{alg:gwrk} using three choices of $\beta \in \{ 0.10, 0.35, 0.50\}$ and explicitly compute the right- and left-hand sides (labeled RHS and LHS in the figures, respectively) of~\eqref{eq:medexp_condition}. Parameters and system construction is as in Figure~\ref{fig:medexp_varybeta}, except for the construction of the corruption vector, which has the same non-zero magnitude across all corruptions (value 1 in the first column and 0.1 in the second column). Figure~\ref{fig:medexp_condition} shows the results of this experiment and shows that, in many iterations,~\eqref {eq:medexp_condition} holds, especially in later iterations and when $\beta < 0.5$.

\subsubsection{Inverse Magnitude Weighting Functions}\label{subsubsec:inverse_type}

We finally focus on the inverse weighting functions.  In this section, we note that $f_t(0) = 0$ since the weighting function is defined with $f_t(r) = 0$ for $r \le 0$, and thus the conclusion of Theorem~\ref{thm:nonincreasing_f_main} does not hold for this class of weighting functions. Despite this, we demonstrate numerically that such weighting functions are still feasible numerically in Section~\ref{sec:experiments} and show that contraction is possible. We prove a single-iteration bound on the expected error of GWRK with inverse weighting functions in Corollary~\ref{cor:inverse-mag-weight}.  We note that this bound is likely quite hard or impossible to iterate due to the dependence on iteration-specific parameters that do not yield easy bounds.  Regardless, this bound may provide contraction in certain instances and may explain the initial convergence we see in numerical experiments with these highly conservative methods.

\begin{corollary}[Inverse convergence]
Suppose $\vb = \mA\vx^* + \vc$ where $\mA \in \mathbb{R}^{m \times n}$ is row-normalized and $C = \text{supp}(\vc)$ satisfies $|C| = \beta m$ and $\beta<1/2$.
Let $\vx^{(t)}$ denote the iterates produced in the $t^{th}$ iteration of Algorithm~\ref{alg:gwrk} where
$f_t(r) = \frac{1}{r^q}$ for $r>0$ with $1 \leq q \leq 2$. Further assume that there is only one zero value in the residual vector, i.e., that $|\{i: r_i^{(t-1)} = 0\} |=1$. Then the iterates satisfy
\begin{align*}
  \mathbb{E}_{t-1}\big[\|\vx^{(t)}-\vx^*\|^2\big] &\leq \left( 1 - (1-\beta)\epsilon_t  \right) \|\vx^{(t-1)}-\vx^*\|^2  \\ 
  & \quad\quad + \frac{\beta}{(1-2\beta)^{q/2}} \frac{m}{(m-1)^{q/2}} \frac{s_{\max}^q}{\delta_t^{q-1}} \| \vx^{(t-1)} - \vx^*\|^q \left(\| \vc\|_\infty +  \| \vx^{(t-1)} - \vx^*\|  \right).
\end{align*}
where $\epsilon_t, \delta_t > 0$ are defined as: 
\begin{align*}
    \epsilon_t&:= \frac{\min_{i \in [m],  |r_i^{(t-1)}| \neq 0} |r_i^{(t-1)}|^2 }{\| \vx^{(t-1)} - \vx^*\|^2},\\ 
    \delta_t &:= \min_{j \in C_F^{(t)}} \left( |c_j| - |\langle \va_j,\vx^{(t-1)} - \vx^* \rangle|]\right),
\end{align*}
and $C_F^{(t)}$ is the `fast' set of corrupted indices defined in~\eqref{eq:fast_C}.
\label{cor:inverse-mag-weight}
\end{corollary}
\begin{proof}
    We begin our proof from \eqref{eq:three_terms} and proceed by bounding the second and third terms. 

    For the second term, we have that: 
    \begin{align*}
        \frac{1}{\sum_{i=1}^m f_t(|r_i^{(t-1)}|)} &\sum_{j\not\in C, r_j \neq 0} f_t(|r_j^{(t-1)}|) |\langle \va_{j}, \vx^{(t-1)} - \vx^*\rangle|^2 \\ 
        & =  \frac{1}{\sum_{i=1}^m f_t(|r_i^{(t-1)}|)} \sum_{j\not\in C, r_j \neq 0} \frac{1}{|r_j^{(t-1)}|^q }  (r_j^{(t-1)})^2 \\ 
         & \geq \frac{(1-\beta) (m-1) r_{\min}^{2-q} }{\sum_{i=1}^m f_t(|r_i^{(t-1)}|)}\\
         & \geq (1-\beta)r_{\min}^2 \\ 
         & \geq (1-\beta) \epsilon_t \|\vx^{(t-1)}-\vx^*\|^2,
    \end{align*}
where the first inequality follows as $\sum_{i=1}^m f_t(|r_i^{(t-1)}|) \leq \frac{m-1}{r_{\min}^q}$, and $r_{\min}$ is the absolute value of the smallest non-zero component of the residual vector $\vr^{(t-1)}$, and the second inequality follows as there exists $\epsilon_t > 0$ such that $r_{\min}^2 \geq \epsilon_t \|\vx^{(t-1)}-\vx^*\|^2.$

    For the last term of \eqref{eq:three_terms},
    \begin{align*}
        \frac{1}{\sum_{i=1}^m f_t(|r_i^{(t-1)}|)} &\sum_{j\in C} f_t(|r_j^{(t-1)}|) [c_{j}^2 - |\langle \va_{j}, \vx^{(t-1)} - \vx^*\rangle|^2] \\ 
        &=  \frac{1}{\sum_{i=1}^m f_t(|r_i^{(t-1)}|)} \sum_{j\in C} \frac{1}{|c_j - \langle \va_j,\vx^{(t-1)} - \vx^* \rangle |^q} [c_{j}^2 - |\langle \va_{j}, \vx^{(t-1)} - \vx^*\rangle|^2] \\ 
         &\leq  \frac{1}{\sum_{i=1}^m f_t(|r_i^{(t-1)}|)} \sum_{j\in C_F^{(t)}} \frac{1}{|c_j - \langle \va_j,\vx^{(t-1)} - \vx^* \rangle |^q} [c_{j}^2 - |\langle \va_{j}, \vx^{(t-1)} - \vx^*\rangle|^2] \\
          &\leq \frac{1}{\sum_{i=1}^m f_t(|r_i^{(t-1)}|)}  \sum_{j\in C_F^{(t)}} \frac{|c_{j}| + |\langle \va_{j}, \vx^{(t-1)} - \vx^*\rangle|}{|c_j - \langle \va_j,\vx^{(t-1)} - \vx^* \rangle|^{q-1} } . 
         \end{align*}

         Suppose that for all $j \in C_F^{(t)}$ there exists a $\delta_t$ such that 
    $$|c_j| - |\langle \va_j,\vx^{(t-1)} - \vx^* \rangle| \geq  \delta_t >
    0, $$
    where $\delta_t$ controls the gap between the corruption value and the residual value. Then, since $q \geq 1$
         \begin{align*}
          \frac{1}{\sum_{i=1}^m f_t(|r_i^{(t-1)}|)}  &\sum_{j\in C_F^{(t)}} \frac{|c_{j}| + |\langle \va_{j}, \vx^{(t-1)} - \vx^*\rangle|}{|c_j - \langle \va_j,\vx^{(t-1)} - \vx^* \rangle|^{q-1} } \\ 
         &\leq \frac{1}{\sum_{i=1}^m f_t(|r_i^{(t-1)}|)} \sum_{j\in C_F^{(t)}} \frac{|c_{j}| + |\langle \va_{j}, \vx^{(t-1)} - \vx^*\rangle|}{\delta_t^{q-1} } \\ 
         & \leq \frac{r_{\min}^q}{\delta_t^{q-1}} \sum_{j\in C_F^{(t)}} |c_{j}| + |\langle \va_{j}, \vx^{(t-1)} - \vx^*\rangle| \\ 
        &  \leq \frac{r_{\min}^q}{\delta_t^{q-1}}  \beta m \left( \|c\|_{\infty} + \| \vx^{(t-1)} - \vx^*\|  \right), 
    \end{align*}
    where we used the fact that $\frac{1}{r_{\min}^q} \leq \sum_{i=1}^{m} f_t(r_i^{(t-1)})$. We have
\begin{align*}
\mathbb{E}_{t-1}\big[\|\vx^{(t)}-\vx^*\|^2\big] 
&\leq \left( 1 - (1-\beta) \epsilon_t \right) \|\vx^{(t-1)}-\vx^*\|^2 \\
& \quad\quad
+ \beta m \frac{r_{\min}^q}{\delta_t^{q-1}} \left(\| c\|_\infty +  \| \vx^{(t-1)} - \vx^*\|  \right).
\end{align*}

By \cite[Lemma 1]{steinerberger2023quantile}, for any $\tau \geq \frac{1}{m-1}$, $|r_{\min}| \leq Q_\tau(\{|r_i^{(t-1)}|\}_{\{i: r_i^{(t-1)} \neq 0\}})\leq \kappa \| \vx^{(t-1)} - \vx^*\|$,
where $\kappa$ is defined in~\eqref{eq:kappa}.  Thus, choosing $\tau = \beta$ and since $\beta<1/2$, we have
$$|r_{\min}|^q \leq \frac{s_{\max}^q}{(m-1)^{q/2}(1-2\beta)^{q/2}} \| \vx^{(t-1)} - \vx^*\|^q.$$
Thus,  we have
\begin{align*}
\mathbb{E}_{t-1}\big[\|\vx^{(t)}-\vx^*\|^2\big] &\leq \left( 1 - (1-\beta) \epsilon_t \right) \|\vx^{(t-1)}-\vx^*\|^2 \\ & \quad\quad+ \frac{\beta}{(1-2\beta)^{q/2}} \frac{m}{(m-1)^{q/2}} \frac{s_{\max}^q}{\delta_t^{q-1} } \| \vx^{(t-1)} - \vx^*\|^q \left(\| \vc\|_\infty +  \| \vx^{(t-1)} - \vx^*\|  \right).
\end{align*}

\end{proof}

\begin{remark} As a consequence of Corollary~\ref{cor:inverse-mag-weight}, we have the following guarantee when $q=1$: 
\begin{align*}
    \mathbb{E}_{t-1}\big[\|\vx^{(t)}-\vx^*\|^2\big] 
    &\leq \left( 1 - (1-\beta) \epsilon_t \right) \|\vx^{(t-1)}-\vx^*\|^2 \\ 
    & \quad\quad
    + \frac{\beta}{(1-2\beta)^{1/2}} \frac{m}{(m-1)^{1/2}} s_{\max} \| \vx^{(t-1)} - \vx^*\| \left(\| \vc\|_\infty +  \| \vx^{(t-1)} - \vx^*\|  \right),
\end{align*}
and when $q=2$:
\begin{align*}
    \mathbb{E}_{t-1}\big[\|\vx^{(t)}-\vx^*\|^2\big] 
    &\leq \left( 1 - (1-\beta) \epsilon_t \right) \|\vx^{(t-1)}-\vx^*\|^2 \\ 
    & \quad\quad
    + \frac{\beta}{(1-2\beta)} \frac{m}{(m-1)} \frac{s_{\max}^2}{\delta_t} \| \vx^{(t-1)} - \vx^*\|^2 \left(\| \vc\|_\infty +  \| \vx^{(t-1)} - \vx^*\|  \right).
\end{align*}
These choices of $q$ are used throughout the remainder of this paper. 
\end{remark}

\section{Numerical Experiments}
\label{sec:experiments}

We next turn to empirical experiments where we explore the behavior of GWRK with a variety of weighting function, and in several corruption regimes.  
Code for all numerical experiments are publicly available at \url{https://github.com/jamiehadd/GWRK}. Experiments are available in both Python and MATLAB. 

\subsection{Experimental Setup}
\label{ssec:setup}
To evaluate the three weighting function families discussed in Section~\ref{subsec:special_cases}, we construct two sparsely corrupted linear systems with different corruption magnitudes.

We generate a matrix $\mA \in \mathbb{R}^{1000 \times 50}$ and a solution vector $\vx^* \in \mathbb{R}^{50}$, both with i.i.d. standard Gaussian entries. The corrupted right-hand side is defined as
$\vb = \mA\vx^* + \vc$, where $\vc \in \mathbb{R}^{1000}$ is a sparse corruption vector with $200$ nonzero entries, corresponding to a corruption ratio of $\beta = 20\%$.
In the \emph{hard system}, the nonzero entries of $\vc$ are drawn i.i.d. from $\mathcal{N}(0,1)$. In this case, the residuals at initialization do not clearly separate corrupted from uncorrupted rows (left panel of Figure~\ref{subfig:two_sys_distribution}).
In the \emph{easy system}, the nonzero entries of $\vc$ are instead drawn i.i.d. from $\mathcal{N}(100,1)$, and the initial residuals of the corrupted and uncorrupted groups are well separated (right panel of Figure~\ref{subfig:two_sys_distribution}).

On both systems, we test several examples for each of the weighting-function family. The details of the example weighting functions we test are as follows:
\begin{itemize}[topsep=4pt, noitemsep]

    \item The threshold weighting function has the form
    $$
    f_t(r)=\mathbf{1}_{\{r\le\theta^{(t-1)}\}}.
    $$
    We test the constant threshold with 
    \begin{equation}\label{ex_func:cons_threshold}
        \theta^{(t-1)} = \theta=Q_{\scriptscriptstyle 0.5}(|\vr^{(t=0)}|),
    \end{equation}
    and the quantile threshold with three different quantile choices 
    \begin{empheq}[
        left={
            \theta^{(t-1)}=Q_{\tau}(|\vr^{(t-1)}|),
            \;\empheqlbrace
        }
    ]{align}
        \tau &= 0.7, \label{ex_func:qrk0.7}\\
        \tau &= 0.8, \label{ex_func:qrk0.8}\\
        \tau &= 0.9. \label{ex_func:qrk0.9}
    \end{empheq}

    \item  The exponential weighting function has the form
    $$
    f_t(r)=\exp\! \Big(-\frac{r^2}{g(|\vr^{(t-1)}|)} \Big).
    $$
    We test four example functions with
    \begin{empheq}[
    left={g^{(t-1)}=\empheqlbrace}
    ]{align}
        & Q_{\scriptscriptstyle 0.5}(|\vr^{(t-1)}|),
          \label{ex_func:exp_med}\\[2pt]
        & Q_{\scriptscriptstyle 0.5}^{2}(|\vr^{(t-1)}|),
          \label{ex_func:exp_med_sq}\\[2pt]
        & \sigma_{\scriptscriptstyle\mathrm{MAD}}(|\vr^{(t-1)}|),
          \label{ex_func:exp_sig}\\[2pt]
        & \sigma_{\scriptscriptstyle\mathrm{MAD}}^{2}(|\vr^{(t-1)}|).
          \label{ex_func:exp_sig_sq}
    \end{empheq}

    \item The inverse weighting function has the form
    $$f_t(r)=r^{-q}\;\;\mathrm{if}\; r>0,\;\mathrm{and}\;0\;\mathrm{otherwise}.$$
    We test three example functions with 
    \begin{empheq}[
    left={q=\empheqlbrace}
    ]{align}
        & 1.0,
          \label{ex_func:inv_q1.0}\\[2pt]
        & 1.2,
          \label{ex_func:inv_q1.2}\\[2pt]
        & 2.0.
          \label{ex_func:inv_q2.0}
    \end{empheq}
    
\end{itemize}
To reduce the influence of any particular random seed on the reported results, we run each algorithm with five different random seeds in every experiment and plot the run with the median number of iterations to convergence.

\subsection{Results and Observations}
For both the hard and easy systems, we test four threshold-type functions (\eqref{ex_func:cons_threshold}, \eqref{ex_func:qrk0.7}, \eqref{ex_func:qrk0.8}, \eqref{ex_func:qrk0.9}), 
four exponential-type functions (\eqref{ex_func:exp_med}, \eqref{ex_func:exp_med_sq}, \eqref{ex_func:exp_sig}, \eqref{ex_func:exp_sig_sq}), 
and three inverse-type functions (\eqref{ex_func:inv_q1.0}, \eqref{ex_func:inv_q1.2}, \eqref{ex_func:inv_q2.0}).

\begin{figure}[htbp]
    
    \centering
    \includegraphics[width=0.98\linewidth]{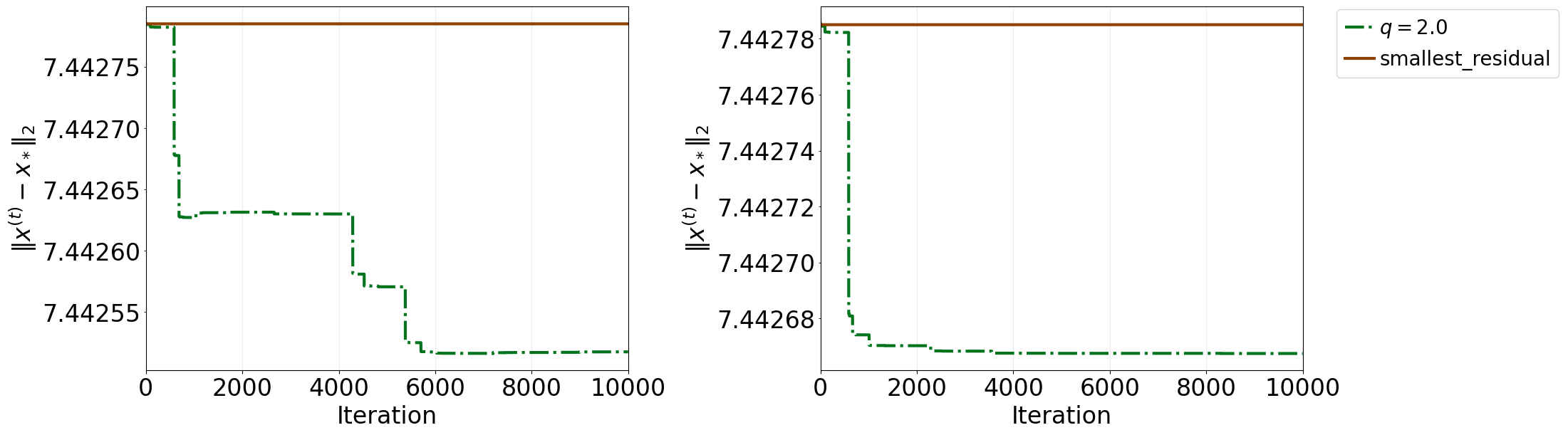}
    \caption{
    Compare the inverse weighting function $f_t(r)=r^{-2}$~\eqref{ex_func:inv_q2.0} against an extremely conservative selection rule: choosing the row which has the smallest nonzero residual, under the hard system (left) and easy system (right).
    System setups are introduced in Section~\ref{ssec:setup}. 
    }
    \label{fig:inverse_vs_smallestresidual}
\end{figure}

\begin{figure}[htbp]
    \centering

    \begin{subfigure}{\linewidth}
        \centering
        \includegraphics[width=0.98\linewidth, height=0.14\textheight, keepaspectratio]{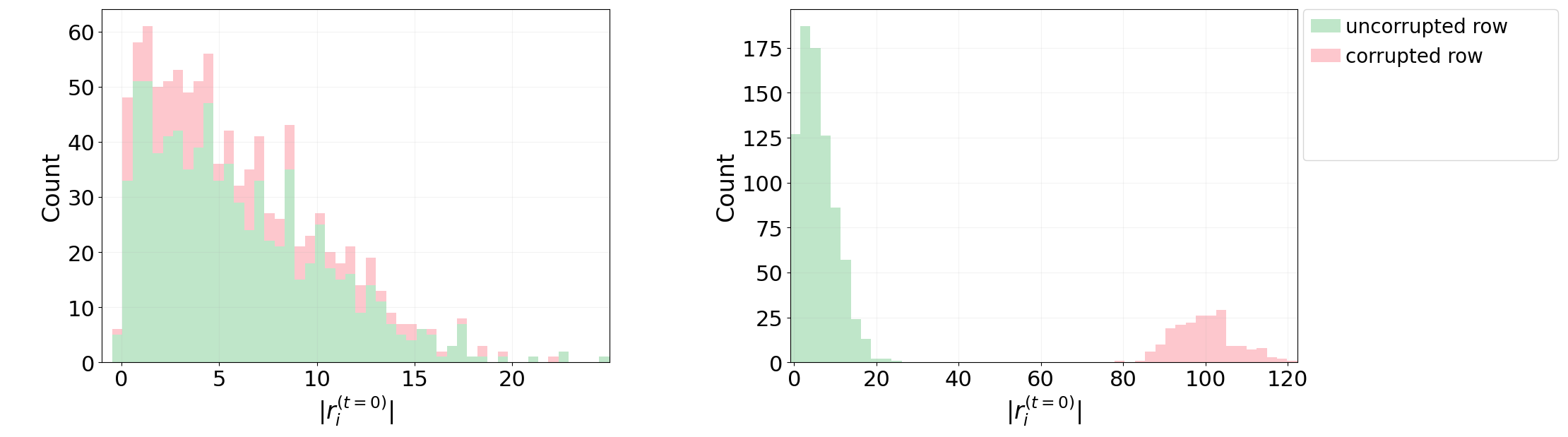}
        
        \caption{
        Distribution of initial residuals: $ \vr^{(0)}=\mA \vx^{(0)} - \vb$.
        }
        \label{subfig:two_sys_distribution}
    \end{subfigure}

    \par

    \begin{subfigure}{\linewidth}
        \centering
        \includegraphics[width=0.98\linewidth, height=0.14\textheight, keepaspectratio]{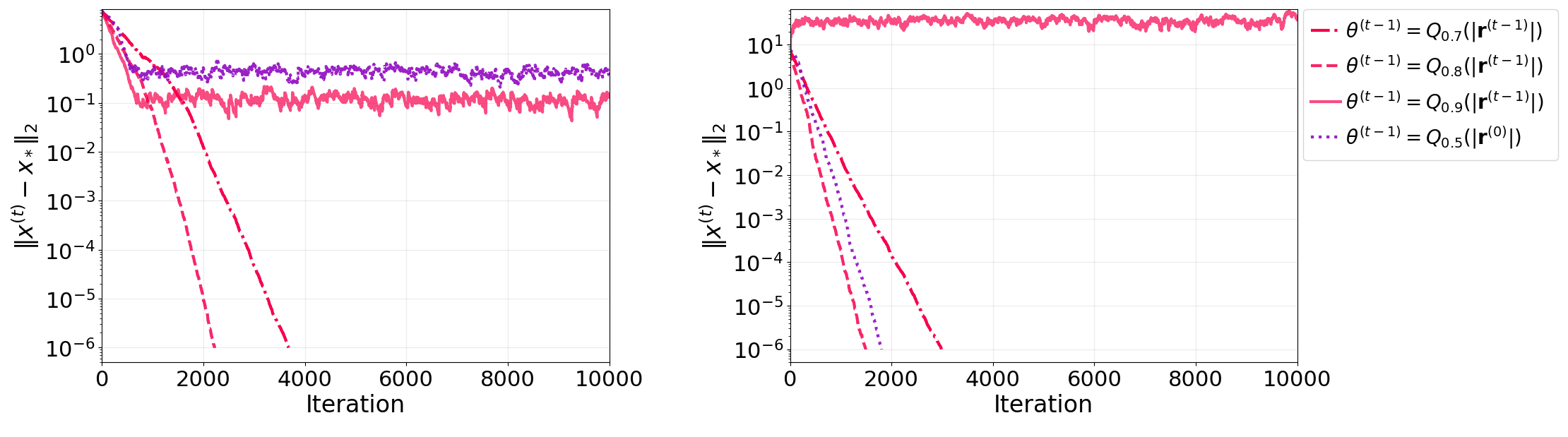}
        
        \caption{
        Threshold weighting:
        $f_t(r)=\mathbf{1}_{\{r\le\theta^{(t-1)}\}}$.
        }
        \label{subfig:two_sys_threshold_type}
    \end{subfigure}

    \par

    \begin{subfigure}{\linewidth}
        \centering
        \includegraphics[width=0.98\linewidth, height=0.14\textheight, keepaspectratio]{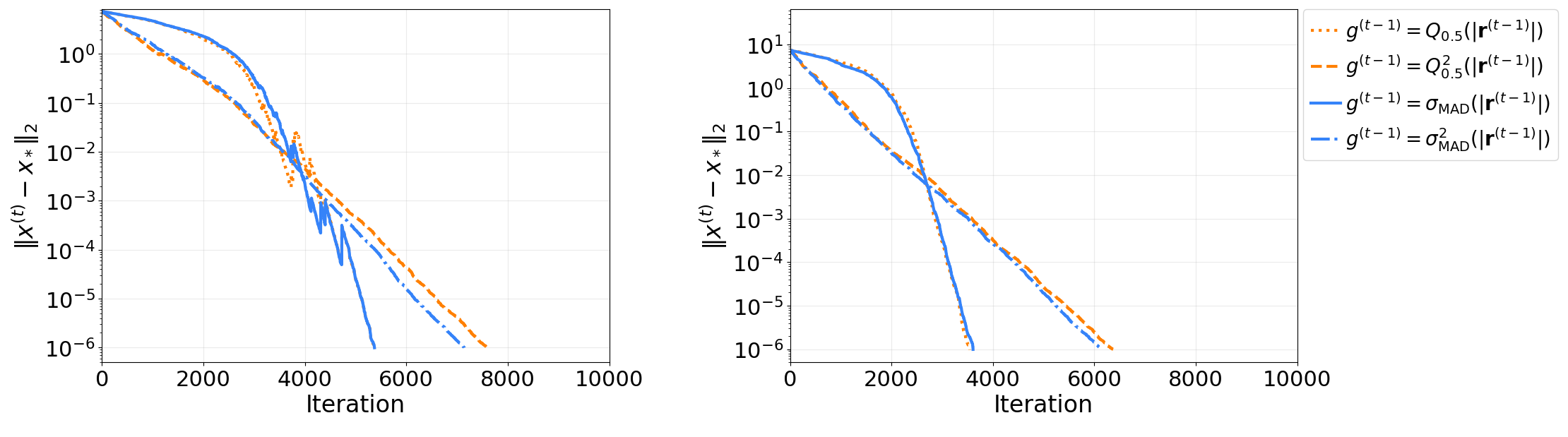}
        
        \caption{
        Exponential weighting:
        $f_t(r)=\exp\!\left(-\frac{r^2}{g(|\vr^{(t-1)}|)}\right)$.
        }
        \label{subfig:two_sys_exp_type}
    \end{subfigure}

    \par

    \begin{subfigure}{\linewidth}
        \centering
        \includegraphics[width=0.98\linewidth, height=0.14\textheight, keepaspectratio]{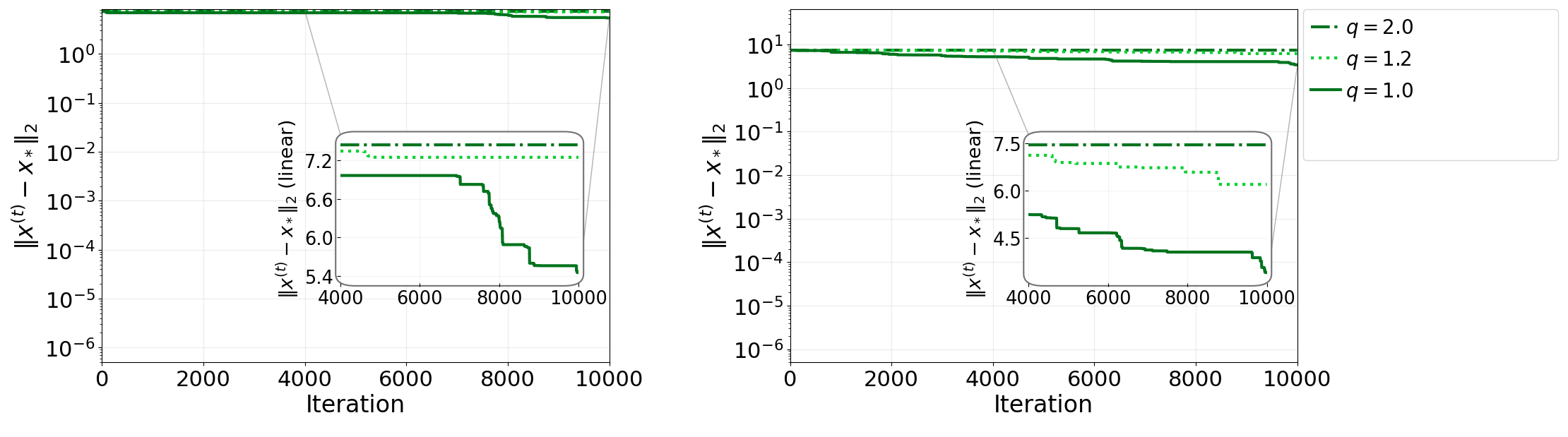}
        
        \caption{
        Inverse weighting:
        $f_t(r)=r^{-q}$ if $r>0$, and $0$ otherwise. 
        }
        \label{subfig:two_sys_inverse_type}
    \end{subfigure}

    \caption{     \textit{Left panel} of each subfigure corresponds to the \textit{hard system}, where the residuals at initialization are unseparable; \textit{right panel} corresponds to the \textit{easy system}, where the residuals at initialization are separable.
    The first subfigure shows the distribution of residuals at initialization.
    The remaining three subfigures present the convergence behavior (approximation error versus iteration) of example functions of the three weighting families listed in Table~\ref{tab:mytable}.
    }
    \label{fig:two_sys_three_type_func}
\end{figure}

\begin{figure}[htbp]
    \centering

    \begin{subfigure}{\linewidth}
        \centering
        \includegraphics[width=0.43\linewidth,             keepaspectratio, height=0.14\textheight, keepaspectratio]{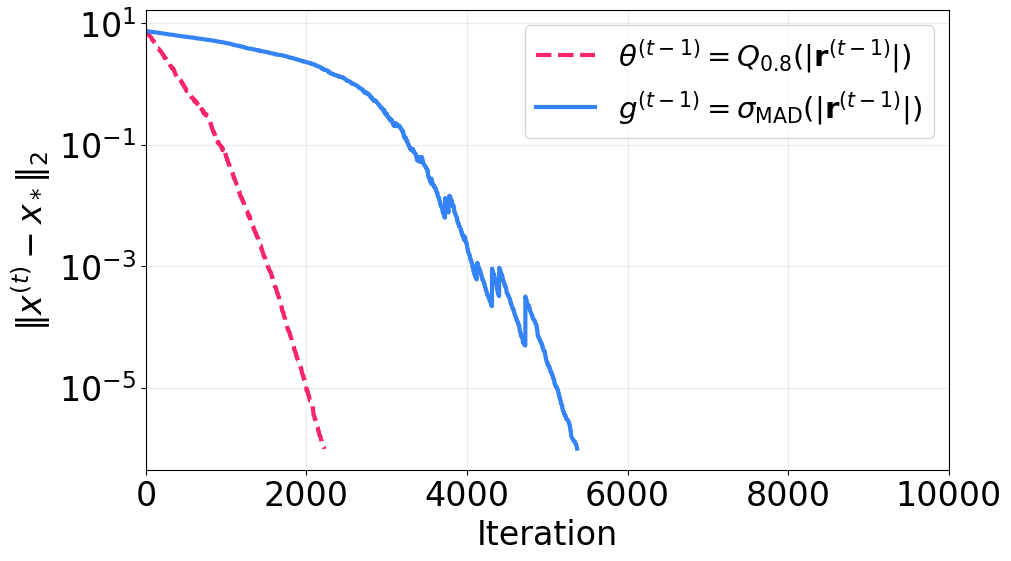}
        
        \caption{
        Two representatives: threshold $\boldsymbol{ Q_{\tau\!=\!0.8}}$ (Fig.~\ref{subfig:two_sys_threshold_type} left) and exponential $\boldsymbol{\sigma_{\scriptscriptstyle \mathrm{MAD}}}$ (Fig.~\ref{subfig:two_sys_exp_type} left).
        }\label{subfig:two_repre_converg}
    \end{subfigure}

    \par 

    \begin{subfigure}{\linewidth}
        \centering
        \includegraphics[width=0.95\linewidth,
            keepaspectratio, height=0.14\textheight, keepaspectratio]{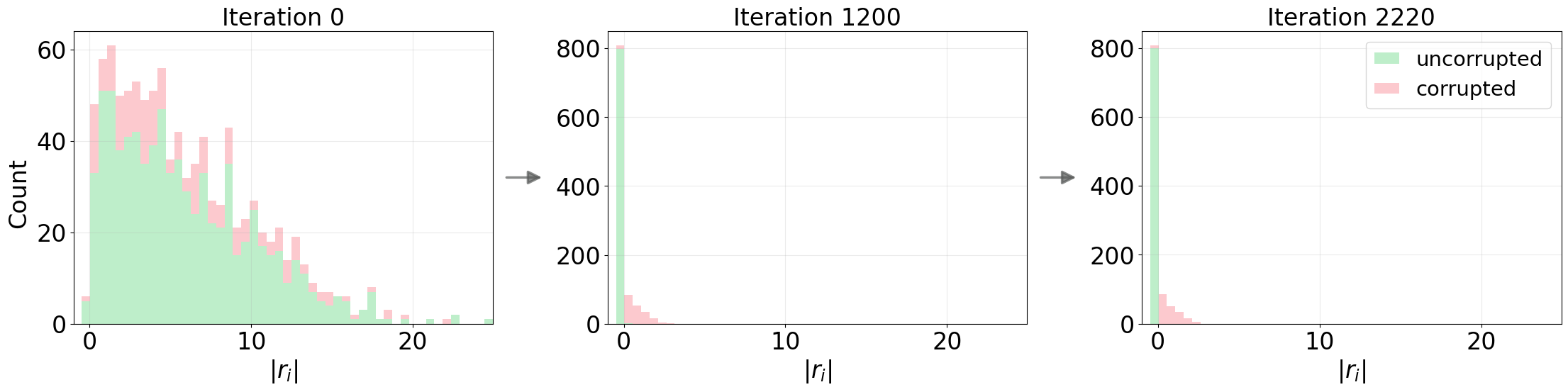}
        
        \caption{
        Residual histograms at the initial/intermediate/final stage for threshold $\boldsymbol{ Q_{\tau\!=\!0.8}}$.
        }\label{subfig:two_repre_qrk_hist}
    \end{subfigure}

    \par 

    \begin{subfigure}{\linewidth}
        \centering
        \includegraphics[width=0.95\linewidth,
            keepaspectratio, height=0.14\textheight, keepaspectratio]{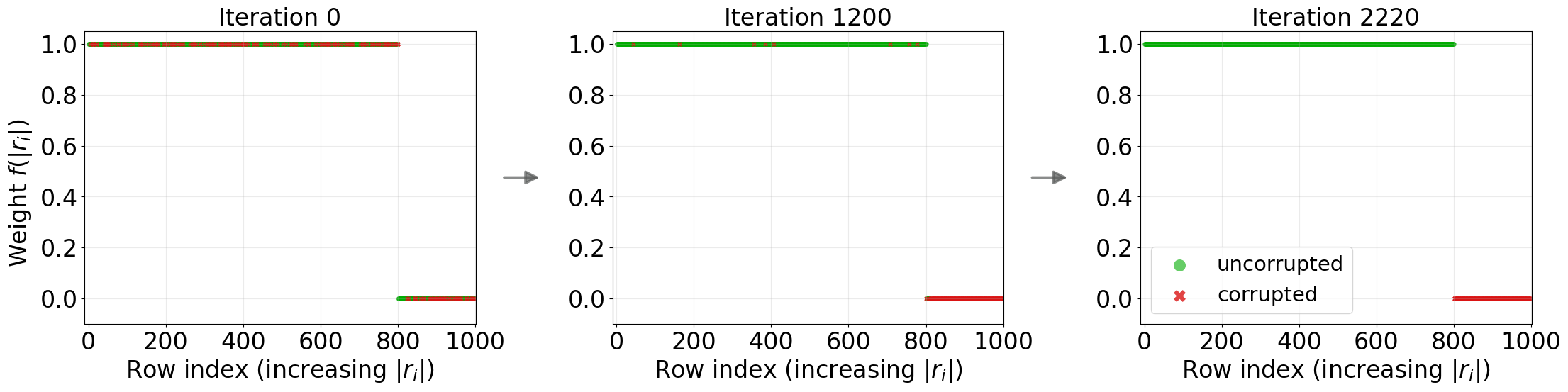}
        
        \caption{
        Weighting dynamics at the initial/intermediate/final stage for threshold $\boldsymbol{ Q_{\tau\!=\!0.8}}$.
        }\label{subfig:two_repre_qrk_weights}
    \end{subfigure}

    \par 

    \begin{subfigure}{\linewidth}
        \centering
        \includegraphics[width=0.95\linewidth,
            keepaspectratio, height=0.14\textheight, keepaspectratio]{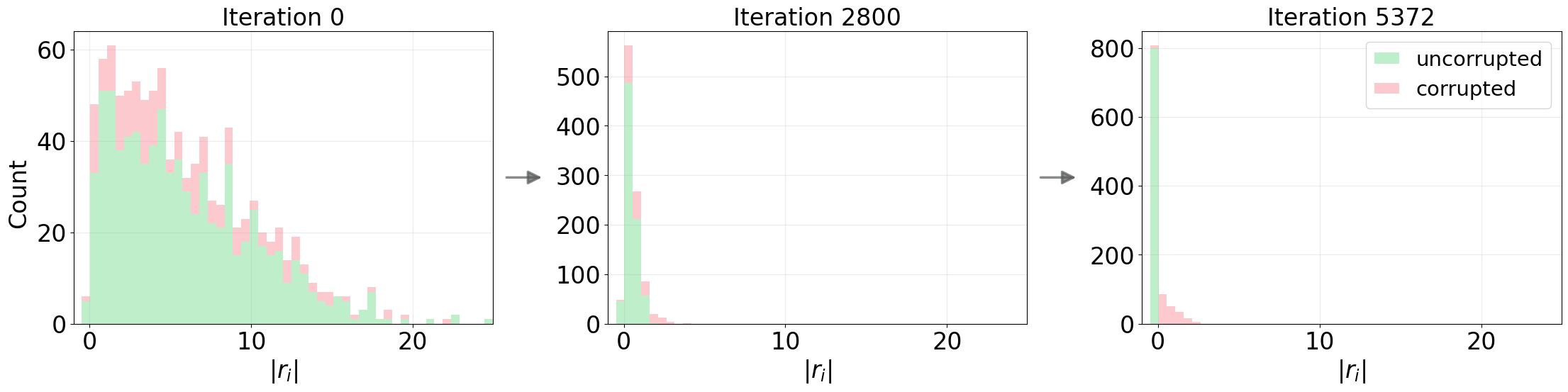}
        
        \caption{
        Residual histograms at the initial/intermediate/final stage for exponential 
        $\boldsymbol{\sigma_{\scriptscriptstyle \mathrm{MAD}}}$.
        }\label{subfig:two_repre_exp_hist}
    \end{subfigure}

    \par 

    \begin{subfigure}{\linewidth}
        \centering
        \includegraphics[width=0.95\linewidth,
            keepaspectratio, height=0.14\textheight, keepaspectratio]{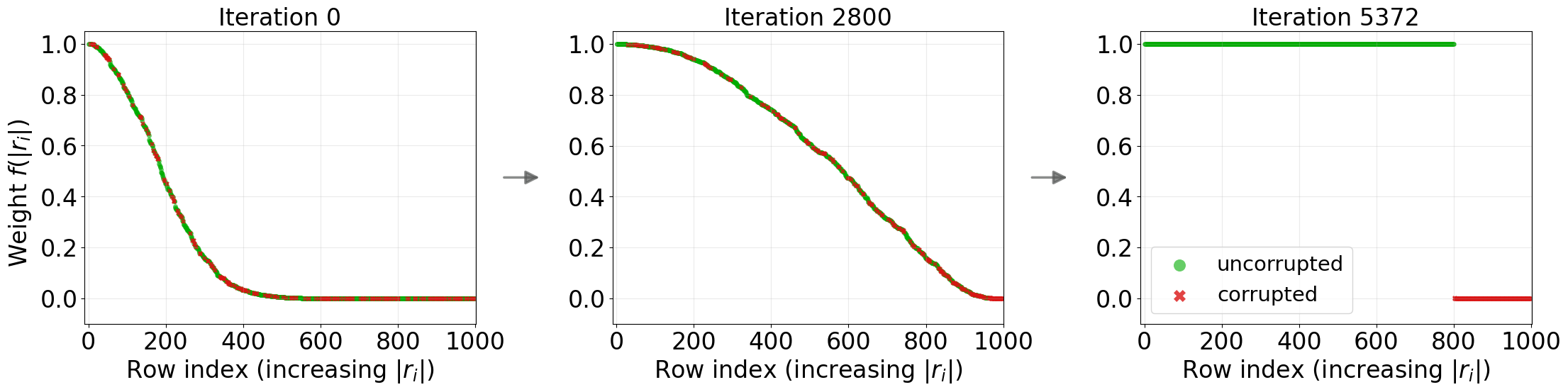}
        
        \caption{
        Weighting dynamics at the initial/intermediate/final phase for exponential 
        $\boldsymbol{\sigma_{\scriptscriptstyle \mathrm{MAD}}}$.
        }\label{subfig:two_repre_exp_weights}
    \end{subfigure}
    
    \caption{In the hard system, evolution of residuals and weights of two representative functions:  threshold-weighting $Q_{\tau\!=\!0.8}$~\eqref{ex_func:qrk0.8} 
    and exponential-weighting $\sigma_{\scriptscriptstyle \mathrm{MAD}}(|\vr^{(t-1)}|)$~\eqref{ex_func:exp_sig}. 
    The exponential weighting exhibits dynamics \emph{similar}
    to those of threshold weighting.
    }
    \label{fig:weighting_dynamics}
\end{figure}

Figures~\ref{fig:inverse_vs_smallestresidual} and~\ref{fig:two_sys_three_type_func} show their convergence behavior, while Figure~\ref{fig:weighting_dynamics} illustrates the evolution of
residuals and weights for two representative functions.
We make following observations based on these results.

\begin{enumerate}
    \item When the threshold is appropriately chosen, \textbf{threshold weighting} converges in the fewest iterations among the three families tested. Both systems have a corruption ratio of $\beta\!=\!20\%$.
    QRK with quantile $\tau\!=\!0.9$ fails to converge (solid pink line in Fig.~\ref{subfig:two_sys_threshold_type}), whereas the choices $\tau\!=\!0.8$ and $\tau\!=\!0.7$ yield rapid convergence (dashed and dash-dotted pink lines, respectively). 
    The function using a constant threshold fixed at $\mathrm{median}(|\vr^{(t=0)}|)$ (dotted purple line in Fig.~\ref{subfig:two_sys_threshold_type}) fails to converge in the hard system, where the initial residuals of corrupted and uncorrupted rows are not separable. In contrast, it performs well in the easy system,  converging slightly slower than QRK with $\tau\!=\!0.8$.

    \item All four evaluated \textbf{exponential weighting} functions exhibit good convergence behavior in both systems, suggesting robustness to
    the initial residual configurations. 
    These functions require no manual parameter tuning, although they generally converge more slowly than the threshold weighting with a suitably chosen threshold.
    Interestingly, we observe that some exponential functions exhibit convergence behavior similar to that of QRK.

    To examine this behavior more closely, we compare two representative functions: QRK with $\tau\!=\!0.8$~\eqref{ex_func:qrk0.8} and exponential $\sigma_{\scriptscriptstyle \mathrm{MAD}}$ weighting~\eqref{ex_func:exp_sig}. Figure~\ref{fig:weighting_dynamics} presents their convergence curves, together with snapshots of the residual distributions and row weights at the initial, intermediate, and final stages.
    Fig.~\ref{subfig:two_repre_converg} shows that although exponential weighting is more conservative and takes more
    iterations to converge, its convergence curve exhibits a slope similar
    to that of QRK during the second half of the iterations.
    The residual and weight dynamics further illustrate these similarities.
    Figs.~\ref{subfig:two_repre_qrk_hist} and~\ref{subfig:two_repre_exp_hist} show that, for both functions, the residual distributions of corrupted and uncorrupted rows overlap
    initially but become separated as the iterations proceed. Differences between the intermediate residual snapshots also reflect the slower convergence of exponential weighting relative to QRK.
    Figs.~\ref{subfig:two_repre_qrk_weights} and~\ref{subfig:two_repre_exp_weights} show the raw weight $f_t(|r_i^{(t-1)}|)$ assigned to each row $i$ at these three stages, with the corresponding iteration labeled above each subfigure.
    Initially, the two functions exhibit distinct weight profiles:
    QRK has a step-like shape, whereas exponential $g\!=\!\sigma_{\scriptscriptstyle \mathrm{MAD}}$ has an S-shaped profile. By the final stage, however, the profiles become similar: uncorrupted rows with small residuals receive high weights, while corrupted rows with large residuals receive zero or nearly zero weights. 
    This is consistent with our discussion on the conservative nature of these functions at the beginning of Section~\ref{subsec:special_cases}.

    \item \textbf{Inverse-weighting} is the most conservative of the three families tested and converges the slowest.
    Larger values of $q$ make the algorithm more conservative and further slow convergence, as illustrated by the zoomed-in insets in
    Fig.~\ref{subfig:two_sys_inverse_type}.
    To further investigate this behavior, we compare inverse weighting with $q\!=\!2.0$ against an extreme selection rule that always chooses the
    row with the smallest nonzero residual.
    The results are shown in Figure~\ref{fig:inverse_vs_smallestresidual}.
    Although convergence is slow, the approximation error of inverse weighting with $q\!=\!2.0$ is indeed decreasing, which validates our analysis in Section~\ref{subsubsec:inverse_type}. 

\end{enumerate}

\subsubsection{Exponential-weighting, Phase Transition}

\label{subsec:numerical_median_exp}

Corollary~\ref{cor:medianexp} presents an upper bound for the expected approximation error when using the median-based exponential downweighting function~\eqref{eq:ff_medexp}. In this experiment, we investigate empirically whether a convergence horizon will be observed in different runs of Algorithm~\ref{alg:gwrk}. Here, we vary the corruption rate $\beta$ while keeping all other parameters fixed. 

The experiments conducted in this subsection use the hard system setup described in Section \ref{ssec:setup} with the exception that instead of running experiments on a fixed corruption ratio $\beta = 0.20$, we vary $\beta = 0.10, 0.30, \dots, 0.90$. Additionally, we run the same experiment on an even finer mesh of corruption ratios where $\beta = 0.50, 0.52, \dots, 0.60$. We represent the family of exponential weighting functions by the exponential squared median function (\ref{ex_func:exp_med_sq}). With this function, we run Algorithm \ref{alg:gwrk} over three random seeds, and plot the run with the median number of iterations.

\begin{figure}[htbp]
    \centering
    \includegraphics[width=\linewidth]{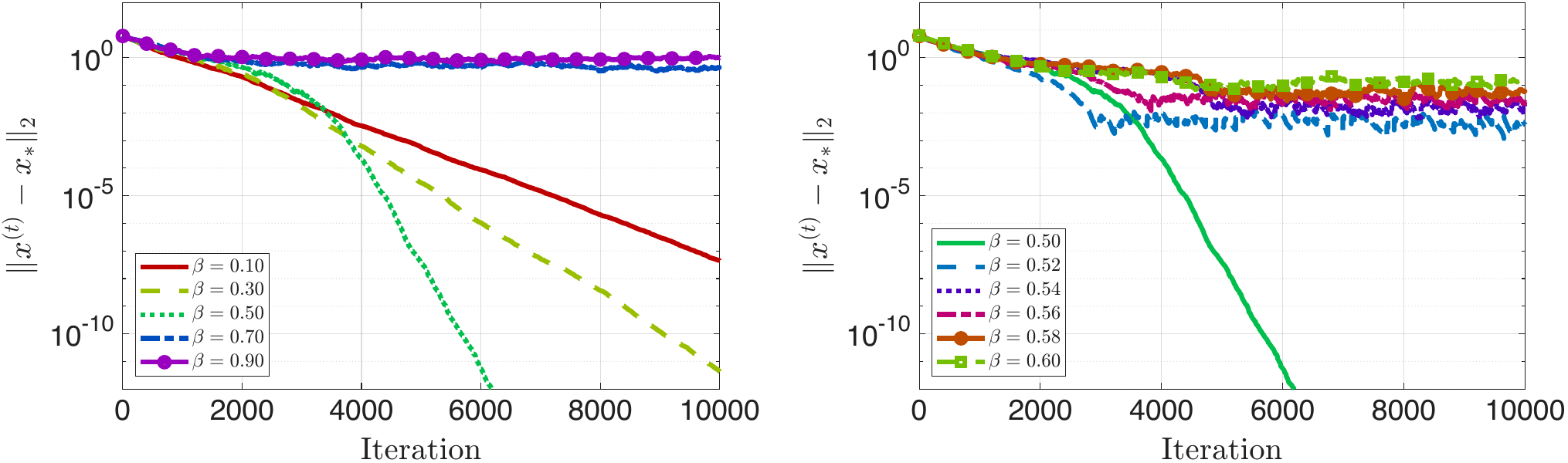}\caption{Performance of Algorithm~\ref{alg:gwrk} with exponential squared median down-weighting (\ref{ex_func:exp_med_sq}) for varying choices of corruption ratio $\beta$.}
    \label{fig:medexp_varybeta}
\end{figure}

Figure~\ref{fig:medexp_varybeta} shows the results of this experiment. The left subplot shows the results for $\beta = 0.10, ..., 0.90$ and the right subplot shows the results for $\beta = 0.50, 0.52, ... 0.60$. The left subplot shows that there is no observed convergence horizon for the corruption rates of interest in Corollary~\ref{cor:medianexp}, $\beta < 0.50$. In addition, the left plot shows that larger corruption rates up to $\beta = 0.50$ actually lead to faster convergence. In particular, we observe that $\beta = 0.50$ achieves the lowest error most quickly. Moreover, the right subplot shows that there is a \emph{phase transition} that occurs when $\beta$ is greater than 0.50. Even increasing $\beta$ slightly (from 0.50 to 0.52) leads to a convergence horizon empirically.

\section{Conclusions and Discussion}

In this work we introduced the generalized weighted randomized Kaczmarz method (GWRK), a framework in which rows are sampled with probability proportional to a user-defined function $f$ of the absolute residuals. In contrast to prior weighted schemes, which favor large magnitude residuals to accelerate convergence, we consider the GWRK sampling scheme to avoid rows whose equations are likely corrupted. Theorem~\ref{thm:nonincreasing_f_main} bounds the expected squared error after one iteration for any non-increasing $f$, separating a contraction term from a convergence-horizon term governed by the weight placed on corrupted rows. We specialized this bound to thresholded, exponential, and inverse-magnitude weighting functions, and the experiments in Section~\ref{sec:experiments} show the resulting methods in practice.

Thresholded weights, which include QRK as a special case, admit the strongest guarantees (Corollary~\ref{cor:thresholded_convergence}) but require a threshold chosen with some knowledge of the corruption rate $\beta$ or magnitudes. Exponential weights scaled by the median or by $\sigma_{\mathrm{MAD}}$ of the residuals are parameter-free and converged in our numerical experiments (Figure~\ref{fig:two_sys_three_type_func}), whereas quantile thresholds with $\tau$ above $1-\beta$ did not. Their weights eventually resemble those of QRK with a well-chosen quantile (Figure~\ref{fig:weighting_dynamics}), but they are more conservative and require more iterations. Inverse-magnitude weights are the most conservative and the slowest.

Our results describe a trade-off between how strongly conservative methods avoid sampling a corrupted equation, and speed of convergence. Roughly stated, the less a weighting function needs to know about the corruption, the more conservative it must be, and the more slowly it converges and the harder it is to prove a horizon-free guarantee.

\subsection{Future Work}
We see several directions for future work.  The first are based upon limitations of our current results.  The contraction condition in Corollary~\ref{cor:thresholded_convergence} depends on worst-case subset singular values and may be restrictive when $\beta$ or the dimension $n$ is large.  Similarly, the bound for inverse-magnitude weights of Corollary~\ref{cor:inverse-mag-weight} is a very restrictive one-step result that we do not expect to iterate.  New analysis techniques would be necessary to improve upon these current limitations.

The bound for GWRK with an exponential weighting function (Corollary~\ref{cor:medianexp}) applies only once the corrupted residuals are separated from the uncorrupted ones, and its contraction factor depends on the iterate-dependent matrix $\tilde{\mathbf{A}}_{C^c}$.  However, the horizon term in Corollary~\ref{cor:medianexp} decays faster than any power of $\|\mathbf{e}^{(t-1)}\|$, so the main obstacles to a horizon-free guarantee are the lack of a uniform lower bound on the contraction factor and the absence of an analysis before the residuals separate. Remark~\ref{rem:med_lower_bound} is a starting point for addressing both.

GWRK was developed for the corrupted linear system case where an otherwise consistent linear system is perturbed by a sparse additive corruption vector.  Natural extensions include systems with both dense noise and sparse corruption, time-varying corruption, and block or subsampled versions of GWRK analogous to those developed for QRK~\cite{cheng2023block,haddock2023subsampled}. Additionally, each iteration requires computing a median or $\sigma_{\mathrm{MAD}}$ of the residual vector, so reducing this cost, for example by incremental updates or by estimating these statistics from a subsample, is important for very large systems.  Finally, an exciting direction of work could adaptively combine weighting techniques.  For instance, one could begin with a conservative parameter-free weight and tighten toward a threshold as the residuals separate, providing an estimate of a quantile level $\tau_t$ from the residual distribution in later iterations.

Another exciting direction for future work is to tighten the relationship between GWRK and robust M-estimation statistical theory.  These approaches use the residuals in different ways. M-estimation minimizes $L_\rho(r_{i})=\sum_{i=1}^m \rho(r_i)$ for a flexible loss $\rho$, which may be solved by iteratively reweighted least squares (IRLS) for suitable choices of $\rho$. GWRK, in contrast, interprets a function of the residuals as a
\emph{sampling distribution} over rows.  Understanding the link between these two approaches could be a fruitful direction of work.  

\section*{Statements and Declarations}

The authors have no relevant financial or non-financial interests to disclose. 

\paragraph{Funding} The authors were supported in part by NSF grant DMS-2520375 that funded the Research Collaboration Workshop in Science of Data and Mathematics (WiSDM) at the University of North Carolina, Chapel Hill, during August 4–8, 2025. 
J.H.\ was partially supported by NSF CAREER \#2440040. A.M.\ was partially supported by Simons Grant MPS-TSM-00014026 and NSF \#2606335. 

\bibliographystyle{plainnat}
\bibliography{references}

\end{document}